\documentclass[11pt]{article}
\usepackage[margin=1in]{geometry}
\usepackage{url}
\usepackage{numprint}
\usepackage{setspace}                
\usepackage{graphicx,circledsteps}
\usepackage[numbers]{natbib}
\usepackage{appendix}
\usepackage{hyperref}
\usepackage{amssymb,amsmath,bbm,float,amsthm}
\usepackage{epstopdf}
\usepackage{xcolor}
\usepackage{enumerate}
\usepackage{prodint}

\usepackage{xr}

\newcommand{\indep}{\perp \!\!\! \perp}

\newtheorem{proposition}{Proposition}

\begin{document}

\title{Fast approximations of pseudo-observations in the context of right-censoring and interval-censoring}
\author{Olivier Bouaziz$^{1}$ }
\date{$^1$Université Paris Cité, CNRS, MAP5, F-75006 Paris, France \\}
\maketitle


\begin{abstract}
In the context of right-censored and interval-censored data we develop asymptotic formulas to compute pseudo-observations for the survival function and the Restricted Mean Survival Time (RMST). Those formulas are based on the original estimators and do not involve computation of the jackknife estimators. For right-censored data, Von Mises expansions of the Kaplan-Meier estimator are used to derive the pseudo-observations. For interval-censored data, a general class of parametric models for the survival function is studied. An asymptotic representation of the pseudo-observations is derived involving the Hessian matrix and the score vector. Theoretical results that justify the use  of pseudo-observations in regression are also derived. The formula is illustrated on the piecewise-constant-hazard model for the RMST. The proposed approximations are extremely accurate, even for small sample sizes, as illustrated on Monte-Carlo simulations and real data. We also study the gain in terms of computation time, as compared to the original jackknife method, which can be substantial for large dataset.\\  


\end{abstract}


\section{Introduction}

In order to study censored data in time to event analysis it is common to model the hazard rate. This allows to correctly take into account censoring in the estimation procedure and provides hazard ratio estimates in the framework of proportional hazard models. However, in some contexts other quantities, that have a more direct interpretation related to the studied problem, might be of interest. 
One example is the Restricted Mean Survival Time (RMST) which is defined as the average survival time up to a fixed point. In that case, it is common (see~\cite{zucker1998restricted},~\cite{chen2001causal},~\cite{zhang2011estimating}) to first model the hazard rate using for instance a Cox model, to derive a survival estimator from the estimated hazard rate and to obtain an estimator of the RMST by integrating out this function. This procedure results in a cumbersome computation where it might be difficult to disentangle the effect of each covariate on the RMST. This is a serious drawback for medical applications and 
there is a need for more direct approaches.  There are several other contexts that are concerned by the difficulty of direct modelling of the quantity of interest. This is typically the case for cumulative incidence functions in a competing risk setting or transition probabilities in a multi-state framework. 

Pseudo-observations have been developed in the seminal work of~\cite{andersen2003generalised} to answer this problem. Those pseudo-observations are constructed using the jackknife method from an estimator of the survival function. 
Theoretical results in~\cite{jacobsen2016note} show that the pseudo-observations can then be used as response variables in a regression model for the quantity of interest, such as the conditional RMST, the cumulative incidence functions in a competing risk setting or the transition probabilities in a multi-state framework. This offers the possibility to directly model the quantity of interest and it is often performed by use of a generalised linear model. 

Another more recent areas of development involving pseudo-observations concerns the study of machine learning methods for time to event analysis. In this context, the problematic is similar: one aims at deriving a complex model, based for instance on neural networks, for quantities of interest such as the survival function (see~\cite{zhao2020deep}),  the cumulative incidence function (see~\cite{sachs2019ensemble},~\cite{ginestet2021deep}) or the RMST (see for instance~\cite{zhao2021deep}). The use of pseudo-observations is then appealing since, once the pseudo-observations are obtained, it is possible to directly use any standard machine learning algorithm by considering those pseudo-observations as (non-censored) response variables. 

Methods based on pseudo-observations are also attractive for interval-censored data. With those data, it is challenging to build a regression model based on semi-parametric methods, for quantities of interest such as the RMST. This is due to the lack of informations induced by interval-censoring. As a matter of fact, even in a nonparametric setting it may be problematic to perform estimation of the survival function. 
In this context, one usually relies on the Turnbull estimator or the convex minorant method which were introduced in~\cite{turnbull76} and~\cite{groeneboom92}, respectively. In~\cite{groeneboom92} it has been proved that these estimators achieve the slow rate of convergence of order $n^{1/3}$ and their distribution is not Gaussian and cannot be explicitly computed. In a regression context, for the estimation of the hazard rate, the Cox model with nonparametric baseline was studied in~\cite{huang1995efficient} but again, the baseline survival function has the $n^{1/3}$ slow rate of convergence and the asymptotic distribution of this estimator could not be derived. As a result, it is common to rely on fully parametric models for modelling quantities such as the survival function or the hazard rate in a regression context. 
In~\cite{lindsey1998study} and~\cite{sun07} a Cox model was studied using parametric baselines such as Weibull or piecewise constant. The methods used to perform estimation are based on maximum likelihood theory where the parametric estimators are derived by maximising the likelihood of the observed data. This allows to recover the classical $\sqrt n$ rate of convergence of the parametric estimators. However, the derivation of the estimators is not explicit, even in the absence of covariates, and rely on a maximisation algorithm such as the Newton-Raphson procedure. In~\cite{bouaziz21} a different approach was proposed based on the EM algorithm by considering the true event times as unobserved variables. This method has the advantage that direct estimators can be computed in the E-step of the algorithm when no covariates are present, which results in a stable and robust estimation procedure. All the aforementioned methods consider estimation of the survival function or of the hazard rate through proportional hazard assumptions, but they are not suited for direct modelling of the RMST, in a regression context. However, this can be achieved by using the pseudo-observations approach. In~\cite{sabathe2020regression}, an illness-death model was considered, and conditional transition probabilities or RMST were computed based on this approach. In order to compute the pseudo-observations, the cumulative transition intensities were estimated using either a penalised spline approach or assuming a Weibull distribution. Similarly, in~\cite{nygaard2020regression} pseudo-observations were computed using a spline approach in order to estimate parameters related to the cumulative incidence function in a competing risk setting.


The key concept about pseudo-observations is that they are built based on the unconditional jackknife estimator of the quantity of interest. While applying the jackknife is straightforward in practice, a limitation of this method comes from the computation burden of calculating the initial estimator $n$ times, where $n$ is the sample size. There exists some R functions designed to improve the computation time, such as the \texttt{jackknife} function in the \texttt{prodlim} package, or the \texttt{pseudo\_independent} function in the \texttt{eventglm} package, which rely on a C++ implementation but the gain is limited as the initial estimator still needs to be implemented $n$ times. The computational burden is particularly important for interval-censored data where there is no direct calculation of the estimators, even in the absence of covariates. In this paper, we develop approximated formulas for pseudo-observations where the jackknife technique does not need to be implemented. In our formulas, the pseudo-observations can be directly computed based on the initial estimator. In the case of right-censored data, we provide formulas based on Von Mises expansion of the Kaplan-Meier estimator. In the case of interval-censored data, we derive general formulas for parametric models that only involve the original estimator, the score function and the Hessian of the density. 
Those formulas are approximations of the original jackknife procedure in the sense that they are equal to the original pseudo-observations up to a remainder term that tends towards $0$ as $n$ tends to infinity. 

However, they turn out to have a very high precision even for moderate sample sizes. Since they only involve the original estimator, the score vector and the Hessian matrix in a parametric context, they are extremely fast to compute thus resulting in a drastic reduction of time.

In the next section, we present a brief summary on the pseudo-values approach. In Section~\ref{sec:RC} we develop asymptotic formulas for computing pseudo-observations of the survival function and the RMST in the context of right-censored data. The case of interval censored data is studied in Section~\ref{sec:IC}. We first discuss the context of nonparametric estimation of the survival function in Section~\ref{sec:npIC}. Then the asymptotic pseudo-observations formulas are developed for general parametric models in Section~\ref{sec:param_model}. In Section~\ref{sec:theo_valid}, theoretical validation of pseudo-observations for parametric models are provided: those results show that the conditional expectation of pseudo-observations approximate the conditional expectation of the response variable of interest. Simulations studies for modelling the conditional RMST in the context of right-censored or interval-censored data are conducted in Section~\ref{sec:simu} where precision and computation time of the approximate formulas are evaluated. Finally, two real data are analysed using the proposed methodology in Section~\ref{sec:realdata}.




\section{Backgrounds on pseudo-regression estimation methods}\label{sec:backgrounds}

Let $T_1^*,\ldots,T_n^*$ be $n$ independent and identically distributed (i.i.d.)  time to event variables of interest, let $\theta$ be a parameter of the form $\theta=\mathbb E[h(T^*_i)]$, where $h$ is a known function. Then introduce $Z_1,\ldots,Z_n$ $n$ i.i.d. covariates and define the conditional expectation $\theta_{(l)}=\mathbb E[h(T^*_l)\mid Z_l]$. We further assume there exists an invertible link function $g$ such that $g(\theta_{(l)})=Z_l^{\top}\beta$, where $\beta$ is a vector of regression parameters of interest. Instead of observing the $T_i^*$'s one usually observes a sample $X_1,\ldots,X_n$ of i.i.d. variables, from which an estimator $\hat\theta$ is constructed. The $l^{\text{th}}$ pseudo-observation is then given by:
\begin{align}\label{eq:pseudo-obs}
\hat\theta_{(l)}=n\hat\theta-(n-1)\hat\theta^{(-l)},
\end{align}
where $\hat\theta^{(-l)}$ is the jackknife estimator of $\hat\theta$, that is the estimator $\hat\theta$ computed on the sample where the $l^{\text{th}}$ observation has been removed.

It has been suggested (see~\cite{andersen2003generalised}) to estimate $\beta$ based on the estimating equation
\begin{align*}
U(\beta)=\sum_{l=1}^n \left(\dot{\theta}_{(l)}\right)^{\top}V_l^{-1}(\hat\theta_{(l)}-\theta_{(l)})=0,
\end{align*}  	
where $\dot{\theta}_{(l)}$ denotes the derivative with respect to $\beta$ of $\theta_{(l)}=g^{-1}(Z_l^T\beta)$ and $V_l$ is a weight matrix. As a result, the estimator $\hat\beta$ verifies the equation $U(\hat\beta)=0$ and it has been suggested to use a sandwich estimator to estimate the variance of $\hat\beta$ (see for instance~\cite{andersen2003generalised} for more details).

In the context of right-censored data, where the $\theta$ parameter is the survival function evaluated at some time point and $\hat \theta$ is its Kaplan-Meier estimator, it has been proved in~\cite{graw2009pseudo} and~\cite{jacobsen2016note} that the resulting estimating function has a mean asymptotically equal to zero. More specifically, one observes $T=\min(T^*,C)$, $\Delta=I(T^*\leq C)$, where $C$ is a right-censoring variable and we define $X=(T,\Delta)$. Let $\theta=S(t)=\mathbb P(T^*>t)$, for $t\in[0,\tau]$. We further assume: 
\begin{itemize}
\item[(i)] $C \indep(T^*,Z)$
\item [(ii)] $\mathbb P(T\geq \tau)>0$.
\end{itemize}
We set $X_i=(T_i,\Delta_i)$, $i=1,\ldots, n$ be $n$ i.i.d replications of ($T,\Delta)$. 
The authors have proved that:
\begin{align}\label{eq:pseudo_theo}
\hat\theta_{(l)}=\theta + \dot{\psi}(X_l)+o_{\mathbb P}(1),
\end{align}
where $\dot{\psi}$ is a first order influence function that verifies $\mathbb E(\dot{\psi}(X_l)\mid Z_l)=\theta_{(l)}-\theta$. 
On the other hand,~\cite{jacobsen2016note} have shown that the sandwich estimator used to estimate the variance of $\hat\beta$ is asymptotically biased. However, the authors have concluded that the difference between their corrected variance estimator and the usual sandwich estimator is of minor importance and as a consequence it is customary to use the sandwich estimator for pseudo-regression.

Once the pseudo-observations have been computed, implementation of the estimating equation along with the sandwich variance estimator can be easily performed from the \texttt{geese} function in the \texttt{geepack} R package. 

In this article, we will present approximate formulas for computing pseudo-observations in the context of right-censoring in Section~\ref{sec:RC} and in the context of interval-censoring in Section~\ref{sec:IC}. Instead of directly using Equation~\eqref{eq:pseudo-obs} to obtain the pseudo-observations, we will present approximated formulas that only involve the estimator $\hat\theta$ computed on the whole sample. In both Sections~\ref{sec:RC} and~\ref{sec:IC}, we will first focus our attention on the problem of modelling $\theta_{(l)}=S(t\mid Z_l)$, the conditional survival function evaluated at time $t$ given the covariate $Z_l$. 
The pseudo-observations can be computed using an estimator of $\theta=S(t)$ the unconditional survival function. A standard link function is $g(\cdot)=\log(-\log(\cdot))$ which gives rise to the Cox model. More complex functions can be chosen for $g$, such as neural-networks (see for instance~\cite{zhao2019dnnsurv},~\cite{zhao2020deep},~\cite{feng2021bdnnsurv}), which will provide performant prediction methods for the conditional survival function. 
Based on those results we will also consider the problem of modelling 
\begin{align}\label{eq:RMST}
\theta_{(l)}=\mathbb E(T^*\wedge \tau\mid Z_l)=\int_0^{\tau}S(t\mid Z_l)dt,
\end{align}
 for some $\tau>0$. This allows to estimate the RMST in a regression context by considering for instance the identity function for $g$ or again more complex link functions such as neural-networks (see for instance~\cite{zhao2021deep}).  In the context of right-censored data only, $\theta$ will be estimated based on the Kaplan-Meier estimator and in the context of interval-censored data, it will be estimated based on parametric models. 


\section{Approximate pseudo-observations for right-censored data}\label{sec:RC}

In this section, we use the same notations as in Section~\ref{sec:backgrounds} for right-censored data. We denote by $\hat S$ the Kaplan-Meier estimator of $S$ and we 
define for $l=1,\ldots,n$, the $l^{\text{th}}$ jackknife estimator $\hat S^{(-l)}$ of $\hat S$ as the estimator constructed when omitting the $l^{\text{th}}$ observation $X_l=(T_l,\Delta_l)$. Introduce $\hat H(t)=\sum_i I(T_i\geq t)/n$ and the observed counting process $N_i(t)=I(T_i\leq t,\Delta_i=1)$, where $I(\cdot)$ is the indicator function. Let
\begin{align*}
\hat\Lambda(t)& = \frac 1n\sum_{i=1}^n\int_0^t\frac{dN_i(u)}{\hat H(u)}
\end{align*}
be the standard Nelson-Aalen estimator (see~\cite{ABGK}) of the cumulative hazard function and define the martingale residuals
\begin{align*}
\hat M_l(t)&=N_l(t)-\int_0^t I(T_l\geq u)d\hat\Lambda(u).
\end{align*}
\begin{proposition}\label{KM}
Under Assumptions (i) and (ii) in Section~\ref{sec:backgrounds} the following results hold:
\begin{align*}
\hat S_{(l)}(t)&=n\hat S(t)-(n-1)\hat S^{(-l)}(t)\\
&=\hat S(t)-\hat S(t)\int _0^t \frac{d\hat M_l(u)}{\hat H(u)}+o_{\mathbb P}(1),
\end{align*}
and
\begin{align*}
\int_0^{\tau}\hat S_{(l)}(t)dt &=\int_0^{\tau}\hat S(t)dt-\int_0^{\tau}\int_u^{\tau}\hat S(t)dt\frac{d\hat M_l(u)}{\hat H(u)}+o_{\mathbb P}(1)\cdot
\end{align*}

\end{proposition}

From those formulas it is clear that the pseudo-observations can be approximated from quantities computed on the original sample. In other words pseudo-observations can be computed without performing the jackknife procedure. This results in a drastic reduction of the computation time of those pseudo-observations as illustrated in the simulation section. We will also see that this approximation is very accurate even for moderate sample sizes. Besides, the main interest for using this formula is for large sample sizes, in particular in a machine learning context where computing pseudo-observations is the first step of the procedure before applying algorithms such as neural networks. In those contexts, the order of the sample size is often in millions. Even though Proposition~\ref{KM} is a direct consequence of the results from~\cite{graw2009pseudo} and~\cite{jacobsen2016note}, a separate proof is provided in the Appendix section.

It should be noted that a different approach for obtaining a fast approximation of pseudo-observations has been considered, based on the infinitesimal jackknife method (see~\cite{jaeckel1972}). It has been implemented in the \texttt{survival} package through the \texttt{pseudo} function. Two versions of the infinitesimal jackknife have been implemented, depending on the method used to estimate the survival function. The \texttt{pseudo} function takes as input a survival function that has either been estimated using the Kaplan-Meier estimator or using the Breslow estimator (the exponential of minus the Nelson-Aalen estimator). Both give very similar formulas as the one proposed in this paper. In particular, all three formulas are asymptotically equivalent. Details about the infinitesimal jackknife and its implementation in the \texttt{survival} package are given in the Supporting Information.

\section{Approximate pseudo-observations for interval-censored data}\label{sec:IC}

In this section, we suppose that instead of directly observing $T^*$ we observe a random interval $[L,R]$, $L\geq 0$ and $L\leq R$, which almost surely contains the event time: $\mathbb P (T^*\in[L,R])=1.$ The right end interval is allowed to take the infinite value such that: 
\begin{itemize}
\item if $0<L<R<\infty$ the data are strictly interval-censored,
\item if $0=L<R<\infty$ the data are left-censored,
\item if $0<L<R=\infty$ the data are right-censored,
\item if $0<L=R<\infty$ the data are exactly observed.
\end{itemize}
Using the notations of Section~\ref{sec:backgrounds} the data then consist of i.i.d. replications $X_i=(L_i,R_i)$, $i=1,\ldots,n$. This situation is often called interval-censoring case $2$ (see~\cite{groeneboom1992information}) when exact observations are not allowed and mixed interval censoring (see~\cite{yu2000consistency}) or partly interval censoring (see~\cite{huang1999asymptotic}) otherwise. In order to derive consistent estimators of the survival function under interval censoring one will usually assume independent censoring in the following way (see for instance~\cite{zhang2005regression}): $\mathbb P(T^*\leq t\mid L=l,R=r)=\mathbb P(T^*\leq t\mid l\leq T^*\leq r)$. This supposes that the variables $(L,R)$ do not convey additional information on the law of $T^*$ apart from assuming $T^*$ to be bracketed by $L$ and $R$. 

\subsection{Comments on the nonparametric case}\label{sec:npIC}

It seems appealing to use the same methodology for interval-censored data as in Section~\ref{sec:RC}. 
A natural nonparametric estimator in the context of interval-censored data is the Turnbull estimator which can be seen as an EM estimator and is consequently rather slow to compute. The gain for avoiding computing $n$ times the Turnbull estimator would therefore be highly significant. 

However, in~\cite{groeneboom1992information} and~\cite{wellner1995interval} it has been showed that this nonparametric maximum likelihood estimator converges at the $n^{1/3}$ or $(n\log(n))^{1/3}$ rates. 
Therefore it will not be possible to derive a relation of the following type:
\begin{align*}
\hat\theta = \theta + \frac 1n \sum_{i=1}^n \dot{\psi}(X_i)+o_{\mathbb P}(n^{-1/2}),
\end{align*}
with $\psi$ verifying $\mathbb E(\dot{\psi}(X_l)\mid Z_l)=\theta_{(l)}-\theta$ such as derived in~\cite{jacobsen2016note} and~\cite{graw2009pseudo}. Because if that would be the case, the convergence rate of $\hat\theta$ would be of the order $n^{1/2}$ due to the central limit theorem. However, this result is crucial to derive Equation~\eqref{eq:pseudo_theo} and assess the validity of the procedure. 

An alternative could be to use results from~\cite{huang1999asymptotic} where it is further assumed that $n_1/n$ tends to a positive constant as $n$ tends to infinity, where $n_1$ is the number of  exact observations. Under this assumption, the authors retrieved a $n^{1/2}$ rate of convergence for the nonparametric maximum likelihood estimator which converges toward a centred gaussian process. However, the covariance function of this process is not explicit and can only be determined as the solution of two integrals. This is caused by the construction of the nonparametric estimator that has no closed form but verifies a self-consistency equation. The asymptotic distribution of the nonparametric estimator was derived using results from infinite dimensional M-estimators in~\cite{van1995efficiency}. The same properties in M-estimators could be used here to derive approximated formulas for the nonparametric survival estimator. However, a careful examination of the proofs in~\cite{huang1999asymptotic} shows that such formulas would lead again to implicit expressions of the pseudo-observations in the same form as the asymptotic limit of the nonparametric survival estimator. Since it does not seem possible to approach those expressions in a straightforward manner we will not pursue this idea. We will focus instead in the next section in modelling the survival function using parametric models.



\subsection{Parametric modelling of the survival function}\label{sec:param_model}

We now assume that the common density function of $T_1^*,\ldots,T_n^*$ depends on $\alpha_0\in\Theta\subset \mathbb R^d$, the true model parameter of dimension $d$. We will denote by $ f^*(t;\alpha_0)$, $\lambda(t;\alpha_0)$, $\Lambda(t;\alpha_0)$ and $S(t;\alpha_0)=\exp(-\Lambda(t;\alpha_0))$ the true density, hazard, cumulative hazard and survival functions of $T^*$, respectively. Instead of directly observing the variables of interest, one usually observes a sample of i.i.d. variables $X_1,\ldots,X_n$ which are assumed to have a common density $f(t;\alpha_0)$. We will use the notations $\nabla \log f(t;\alpha_0)$ and $\nabla^2 \log f(t;\alpha_0)$ to represent the score vector and the Hessian matrix of this log-density where the derivatives are taken with respect to the model parameter $\alpha$ and are evaluated at $\alpha=\alpha_0$. The same notations will be used for $f^*(t;\alpha_0)$.
It is important to emphasise the distinction between the notation $f^*$, which represents the density of the true data, and the notation $f$ which represents the density of the observed data. As an illustration, 
in the general context of mixed interval-censored data, $X_i=(L_i,R_i)$ with $0\leq L_i<R_i\leq \infty$ and we have (see~\cite{sun07} or~\cite{bouaziz21})
\begin{align*}
f(X_i;\alpha)&=(S(L_i;\alpha)-S(R_i;\alpha))I(L_i\neq R_i)+\left(\lambda(L_i;\alpha)S(L_i;\alpha)\right)I(L_i= R_i),
\end{align*}
with the slight abuse of notation $S(R_i;\alpha)=0$ if $R_i=\infty$. In the following, we will consider maximum likelihood estimation for the parameter $\alpha_0$ based on the observed variables $X_1,\ldots,X_n$. The results derived in this section are not limited to the case of interval-censored data and can be applied to any parametric framework for incomplete data.

The maximum likelihood estimator $\hat\alpha$ of $\alpha_0$ maximises with respect to $\alpha$ the log-likelihood $\sum_i \log f(X_i;\alpha)$ and, subject to regularity conditions, verifies the following equality:
\begin{align*}
\sqrt n (\hat \alpha-\alpha_0) & = \frac 1{\sqrt n}\left(-\frac 1n \sum_{i=1}^n \nabla^2 \log f(X_i;\tilde{\alpha})\right)^{-1} \sum_{i=1}^n \nabla \log f(X_i;\alpha_0),
\end{align*}
where $\tilde{\alpha}$ lies between $\hat\alpha$ and $\alpha_0$. Let $I=-\mathbb E\left(\nabla^2 \log f(X;\alpha_0)\right)$ be the Fisher information and consider the jackknife version $\hat \alpha^{(-l)}$ of the maximum likelihood estimator. It is then straightforward to write:
\begin{align*}
\sqrt n (\hat \alpha-\alpha_0) & =\frac 1{\sqrt n} \; I^{-1} \sum_{i=1}^n \nabla \log f(X_i;\alpha_0)+\varepsilon_n,\\
\sqrt{n-1} (\hat \alpha^{(-l)}-\alpha_0) & =\frac 1{\sqrt{ n-1}} \; I^{-1} \sum_{i\neq l}^n \nabla \log f(X_i;\alpha_0)+\varepsilon^{(-l)}_n,
\end{align*}
where
\begin{align*}
\varepsilon_n=\frac 1{\sqrt n}\left(\Bigg(-\frac 1n \sum_{i=1}^n \nabla^2 \log f(X_i;\tilde{\alpha})\Bigg)^{-1}-I^{-1}\right)\sum_{i=1}^n \nabla \log f(X_i;\alpha_0),
\end{align*}
and $\varepsilon^{(-l)}_n$ is the jackknife version of $\varepsilon_n$. As a result, the $l^{\text{th}}$ pseudo-observation of $\hat\alpha$ verifies the relation:
\begin{align}\label{eq:pseudo_param}
n\hat\alpha-(n-1) \hat\alpha^{(-l)} = \alpha_0+I^{-1}\nabla \log f(X_l;\alpha_0)+\sqrt n\, \varepsilon_n-\sqrt{n-1}\, \varepsilon^{(-l)}_n.
\end{align}
In the Appendix section, it is proved that the term $\sqrt n\, \varepsilon_n-\sqrt{n-1}\, \varepsilon^{(-l)}_n$ tends towards $0$ in probability as $n$ tends to infinity. This entails that asymptotically, the pseudo-observation of $\hat\alpha$ only depends on the true parameter, the Fisher information and the score vector, this latter quantity being only evaluated at the observation $l$. 
Since $\hat\alpha$ is a consistent estimator of $\alpha_0$ and $\hat I=-\sum_{i=1}^n \nabla^2 \log f(X_i;\hat{\alpha})/n$ is a consistent estimator of $I$, a natural asymptotic approximation for the pseudo-observation of $\hat\alpha$ is simply:
\begin{align*}
\hat\alpha+{\hat I}^{-1}\nabla \log f(X_l;\hat\alpha).
\end{align*}
While this result is interesting on its own, more work needs to be done in order to derive the pseudo-observations of $S(t;\hat\alpha)$. The following proposition is derived based on this latter expression of the approximate pseudo-observation for $\hat\alpha$. The notation $\cdot^{\top}$ is used to denote the transpose of a vector or a matrix.
\begin{proposition}\label{paramsurv}
Under standard regularity conditions for maximum likelihood theory, the following relations hold:
\begin{align*}
n S(t;\hat\alpha) - (n-1) S(t;\hat\alpha^{(-l)}) &=S(t;\hat\alpha) - S(t;\hat\alpha)\nabla \Lambda (t;\hat\alpha)^{\top}{\hat I}^{-1} \nabla \log f(X_l;\hat\alpha)+o_{\mathbb P}(1),
\end{align*}
and for $\tau>0$,
\begin{align*}
&n \int _0^{\tau}\!S(t;\hat\alpha)dt - (n-1) \int _0^{\tau}\!S(t;\hat\alpha^{(-l)})dt\\
 &\quad=\int _0^{\tau}\!S(t;\hat\alpha)dt - \int _0^{\tau}\!S(t;\hat\alpha)\nabla \Lambda (t;\hat\alpha)^{\top}dt\,{\hat I}^{-1} \nabla \log f(X_l;\hat\alpha)+o_{\mathbb P}(1).
\end{align*}
The same results also hold when replacing $\hat\alpha$, $\hat I$ by $\alpha_0$, $I$ in the right-hand side of the equations.
\end{proposition}
The proof of this proposition can be found in the Appendix section. As in Section~\ref{sec:RC}, the main interest in this result lies in the fact that the approximated version of the pseudo-observation only depends on the parameter estimator $\hat\alpha$ and not on its jackknife version. This means that pseudo-observations in parametric models can be obtained without actually computing the $n$ jackknife estimators. Only the estimator of $\alpha_0$, along the Hessian matrix, the gradient of $\Lambda$ and of the log-density are needed. This is particularly interesting in the context of interval-censored data since parametric estimators cannot be derived explicitly and numeric methods must be implemented. 
Two different strategies exist for those types of data: either a direct maximisation of the likelihood can be performed using the Newton-Raphson algorithm (see~\cite{lindsey1998study} and~\cite{sun07} for instance) or the complete likelihood (based on the unobserved true times) can be used through the EM algorithm in order to maximise the likelihood (see~\cite{bouaziz21}). But in either case the method is iterative. 
Also, it should be noted that the Newton-Raphson algorithm requires to compute the score vector and Hessian matrix. 
Therefore the computational cost for implementing the approximated pseudo-observations is similar to the cost of simply computing the pointwise estimate from the Newton-Raphson algorithm. 

Those approximated formulas are general and work for any parametric model. As an illustration, the piecewise-constant hazard (pch) model will be used in the simulation section. This model assumes that the hazard function verifies $\lambda(t;\alpha)=\sum_{k=1}^K \alpha_k I_k(t)$ where $I_k(t)=I(c_{k-1}<t\leq c_k)$, $c_0=0<c_1<\cdots <c_K=+\infty$ represent $K+1$ cuts and $I(\cdot)$ denotes the indicator function. We do not specify precisely the regularity conditions for maximum likelihood theory to hold. However, two important assumptions are first to assume the model identifiable and second to impose the Fisher information to be positive definite in a neighbourhood of the true parameter. For the pch model in the context of interval-censored data, two necessary conditions for those regularity assumptions to hold are:
\begin{align}\label{eq:condIC}
\mathbb P(R<+\infty, [L,R] \cap (c_{k-1},c_k]\neq \emptyset) >0,&\,\forall k=1,\ldots,K,\nonumber\\
\mathbb P(L>c_{k-1})>0,& \,\forall k=1,\ldots,K.
\end{align}
The first assumption is quite natural: in order to estimate $\alpha_k$, the probability that an interval intersects $[c_{k-1},c_k]$ should be positive. The second assumption is necessary for the existence of a maximum of the likelihood function. 
It should be noted that those conditions are also valid when exact observations $L=R$ are allowed. Exact expressions of the score vector and Hessian matrix for the pch model along with the derivation of condition~\eqref{eq:condIC} are detailed in Section~\ref{sec:details_pseudoIC} of the Appendix. Details on the implementation of Proposition~\ref{paramsurv} for the pch model are given in Section~\ref{sec:scoreHessian_pch} of the Appendix.

Precision and computational cost of the approximation for the RMST are evaluated and compared to the actual jackknife version of the pseudo-observations in the simulation section. In particular, it is seen that the approximation is much faster than the jackknife method and is very accurate even for small sample sizes.


\subsection{Theoretical validation of pseudo-observations for parametric models}\label{sec:theo_valid}

In this section we want to investigate if the approximated formula derived in Proposition~\ref{paramsurv} provides valid observations for performing pseudo-regression, similarly to the Kaplan-Meier estimator in the context of right-censored data (see Equation~\eqref{eq:pseudo_theo}). 
In other words, if we set
\begin{align}\label{eq:varphi}
\varphi(X_l;\alpha_0)=\int _0^{\tau}\!S(t;\alpha_0)dt - \int _0^{\tau}\!S(t;\alpha_0)\nabla \Lambda (t;\alpha_0)^{\top}dt\,{I}^{-1} \nabla \log f(X_l;\alpha_0),
\end{align}
we want to investigate under which conditions we may have 
\begin{align}\label{eq:pseudo_param_theo}
\mathbb E(\varphi(X_l;\alpha_0)\mid Z_l)=\mathbb E(T_l^*\wedge\tau\mid Z_l).
\end{align} 
It is easily seen that this equality will generally not hold by considering the simple scenario of exact observations. In that case, $X_i=T_i^*$ and $f=f^*$ is simply the density of the true variable $T^*$. If we further assume for instance that $T^*$ follows a Weibull distribution, with shape parameter $a>0$ and scale parameter $b>0$ (the true parameters are noted $a_0$, $b_0$) such that:
\begin{align*}
f(X;\alpha)=f(T^*;\alpha)=\frac ab \left(\frac {T^*}b\right)^{a-1}\exp\left(-\left(\frac {T^*}b\right)^a\right),
\end{align*}
with $\alpha=(a\;\; b)^{\top}$, then $\nabla \log f(X_l;\alpha)$ will depend on $\log (T^*_l)$, $(T^*_l)^{a-1}$ and $(T^*)^a$, when $a\neq 1$. As a result, $\mathbb E(\varphi(X_l;\alpha_0)\mid Z_l)$ will be a function of $\mathbb E(\log (T^*_l)\mid Z_l)$, $\mathbb E((T^*_l)^{a_0-1}\mid Z_l)$ and $\mathbb E((T^*_l)^{a_0}\mid Z_l)$ when $a_0\neq1$. When $a_0=1$ (the exponential model), then $\mathbb E(\varphi(X_l;\alpha_0)\mid Z_l)$ will be a function of $\mathbb E(T_l^*\mid Z_l)$, but it will still not verify Equation~\eqref{eq:pseudo_param_theo} unless $\tau=\infty$. Performing the same calculation for other distributions, we can similarly conclude that Equation~\eqref{eq:pseudo_param_theo} will not hold in general.

Nevertheless, even though Equality~\eqref{eq:pseudo_param_theo} is not verified in most cases, it is still possible to prove that the formula given by Proposition~\ref{paramsurv} provides a good approximation. The key is to assume that there exists a value $\alpha_z$ of the parameter $\alpha$ such that the conditional distribution of $T^*$ given $Z=z$ follows a distribution with density $f^*(t;\alpha_z)$. In that case, $\mathbb E(T_l^*\wedge\tau\mid Z_l=z)=\int_0^{\tau} S(t;\alpha_z)dt$ and by expanding $S(t;\alpha_z)$ around $S(t;\alpha_0)$ from a Taylor development, we can prove that $\mathbb E(T_l^*\wedge\tau\mid Z_l=z)$ is equal to $\mathbb E(\varphi(X_l;\alpha_0)\mid Z_l=z)$ up to two remainder terms. Those remainder terms measure the distance between $\alpha_0$ and $\alpha_z$, and between the inverse of the Fisher information $I^{-1}$ and the quantity $(-\mathbb E(\nabla^2\log f(X;\alpha)\mid Z=z))^{-1}$ for any $\alpha$ that is on the real line between $\alpha_0$ and $\alpha_z$. Define $\Theta_{z}=\{\alpha\in\Theta: \|\alpha-\alpha_0\|\leq \|\alpha_0-\alpha_z\|\}$ which represents the set of parameters that are on the real line between $\alpha_0$ and $\alpha_z$. In the next proposition, we denote for $k=1,\ldots,K$, by $\alpha_k$, $\alpha_{0,k}$, $\alpha_{z,k}$ the $k$th component of $\alpha$, $\alpha_0$ and $\alpha_z$, respectively.

 \begin{proposition}\label{paramsurv_approx_property}
Assume there exists $\alpha_z$ such that the conditional distribution of $T^*$ given $Z=z$ follows a distribution with density $f^*(t;\alpha_z)$. Assume also there exists $M_z<+\infty$ such that
\begin{align*}
\forall k,k'\in \{1,\ldots,d\},\forall \alpha\in\Theta_z, \frac{\partial^2}{\partial \alpha_k\partial\alpha_{k'}}\int_0^{\tau}S(t;\alpha)dt\leq M_z.
\end{align*}
Then,
\begin{align*}
\mathbb E(T_l^*\wedge\tau\mid Z_l=z)=\mathbb E(\varphi(X_l;\alpha_0)\mid Z_l=z)+R_{1,z}+R_{2,z},
\end{align*} 
where
\begin{align*}
R_{1,z} & \leq \frac {M_z}2 \left(\sum_{k=1}^d \left(\alpha_{z,k}-\alpha_{0,k}\right)\right)^2,\\
R_{2,z} & \leq \max_{\alpha\in\Theta_z}\left|\int_0^{\tau}(\nabla S(t;\alpha_0))^{\top}dt\left(I_{\alpha,z}^{-1}-I^{-1}\right)\mathbb E(\nabla\log(f(X;\alpha_0))\mid Z=z))\right|,
\end{align*}
with
\begin{align*}
I_{\alpha,z}=-\mathbb E(\nabla^2\log(f(X;\alpha))\mid Z=z).
\end{align*}
\end{proposition}
The proposition makes the strong assumption that the conditional distribution of $T^*$ given $Z=z$ follows a distribution with density $f^*(t;\alpha_z)$. While this will not be true in general, it seems reasonable to assume that, if the chosen parametric distribution for $T^*$ is rich enough, there will exist a value of the parameter $\alpha$ such that the parametric distribution is not too far from the conditional distribution of $T^*$ given $Z=z$. This advocates for flexible parametric models such as the pch model or a spline approach such as proposed in~\cite{nygaard2020regression}. 
However, this needs to be imposed for all possible values of $Z$, which again seems reasonable if discrete covariates are considered and the number of those covariates is not too large. 

Besides, it is difficult to evaluate, from the proposition, how large the remainder terms are. Previous experiments on the Weibull distribution in~\cite{sabathe2020regression} using the jackknife approach suggest that the approximation is quite accurate in practice. When using the pch model, we can establish a different type of theoretical result. In the context or right-censored data, we show in the next proposition that if the number of cuts in the pch model tends to infinity, then we exactly retrieve Equality~\eqref{eq:pseudo_param_theo}. We were not able to prove this result in the context of interval-censored data, we conjecture however that this result still holds in this case.



\begin{proposition}\label{paramsurv_pch_property}
Under the context of right-censored data,  if $T^*$ follows the pch model with cuts $c_0=0<c_1<\cdots <c_K=+\infty$ and if we assume standard regularity conditions for maximum likelihood theory then the function $\varphi(X_l;\alpha_0)$ defined in Equation~\eqref{eq:varphi} converges, as $K$ tends to infinity and $max_{|c_{k+1}-c_k|}$ tends to $0$, towards a function $\varphi^{\infty}(X_l;\alpha_0)$ that verifies
\begin{align*}
\mathbb E(\varphi^{\infty}(X_l;\alpha_0)\mid Z_l)=\mathbb E(T_l^*\wedge\tau\mid Z_l).
\end{align*}
\end{proposition}
This result is interesting as it shows that it is theoretically possible for Equality~\eqref{eq:pseudo_param_theo} to hold true when using the pch model. Of course, in practice one has to choose a finite number of cuts. However, there are some strategies to choose the number of cuts from the data. In particular, in~\cite{bouaziz21} the authors have developed a penalised method based on the adaptive-ridge to choose the number of cuts in an efficient way. In the simulation study (Section~\ref{sec:simuIC}), we show that the approximation formula or the original jackknife method provide very similar and very performant results in pseudo-regression, when modelling the distribution of $T^*$ with the pch model. The proofs of Propositions~\ref{paramsurv_approx_property} and~\ref{paramsurv_pch_property} are deferred to the Appendix section.

\section{Simulation studies for the Restricted Mean Survival Time}\label{sec:simu}

We study two different simulation scenarios for the RMST: one with right-censored data and another one with interval-censored data. In the first scenario, the approximate pseudo observations are based on the Kaplan-Meier estimator (using Proposition~\ref{KM}) while in the second scenario they are based on the pch model (using Proposition~\ref{paramsurv}). 
In both settings, the performance of the estimators derived from the approximated formulas and the ones obtained from the standard jackknife method is compared based on $500$ replications. Implementation of the generalised estimation equation is performed through the \texttt{geese} function in the \texttt{geepack} R package. Computation times of the estimators were evaluated on $100$ replications, from $10$ different samples with $10$ replications on each sample using the \texttt{microbenchmark} R package.

\subsection{Right-censored data}

The simulation setting is based on the one in~\cite{wang2018modeling}. We assume that

\begin{align*}
T^*_i=\tilde\beta_0^{\top} Z_i+\varepsilon_i,\quad i=1,\ldots,n,
\end{align*}
with $\tilde\beta_0=(5.5,0.25,0.25)^{\top}$, $Z_i=(1,Z_{i,1},Z_{i,2})^{\top}$, $Z_{i,1}$ and $Z_{i,2}$ are Bernoulli variables with parameter $0.5$ and $\varepsilon_i\sim \mathcal U[-\sigma,\sigma]$, with $\sigma=3$. Under this model it can easily be seen that
\begin{align}\label{simuRC}
\mathbb E(T^*_i\wedge \tau\mid Z_i)=\beta_{00}+\beta_{01} Z_{i,1}(1-Z_{i,2})+\beta_{10} Z_{i,2}(1-Z_{i,1})+\beta_{11} Z_{i,1}Z_{i,2},
\end{align}
where $\beta_0=(\beta_{00},\beta_{01},\beta_{10},\beta_{11})^{\top}$ can be determined computationally using Monte-Carlo samples with size $10$ million. We further set $\tau=6$ which corresponds to the $54.2\%$ quantile of $T^*$ and to the value $\beta_0=(4.98,0.14,0.14,0.27)^{\top}$. Right-censored data were simulated from an exponential distribution with parameter $0.07$ yielding $33\%$ of censoring on average. The results are presented in Table~\ref{simu:RC}.

It is seen that the approximated formula gives similar results as compared to the standard jackknife method for $n=100$. For larger sample sizes, the results are almost identical. We also compared the difference between the two estimators of $\beta_0$ by looking at the standard deviation for all four components and taking the maximum: the maximum value over all four components is equal to $7.06\times 10^{-3}$, $5.90 \times10^{-4}$, $2.16 \times10^{-4}$, $6.88 \times10^{-6}$ for $n=100$, $n=500$, $n=1,000$, $n=10,000$ respectively. This shows that there is very little variations between the estimator computed from the jackknife and the one computed from the approximated formula. In terms of computation times, there is a clear advantage for the approximated formula which goes $14.3$, $27.5$, $25.1$ and $18.7$ times faster for $n=100$, $n=500$, $n=1,000$ and $n=10, 000$ respectively. Clearly the computation time for the original jackknife method is not a linear function of the sample size and the gain for using the approximated method is considerable for large sample sizes. It should be noted that the computation time was evaluated for the pseudo-regression procedure, but it does not include the computation of the initial survival estimator, it only takes into account the computation of the pseudo-observations along with the implementation of the generalised estimating equations. Finally, the infinitesimal jackknife implemented in the \texttt{pseudo} function of the \texttt{survival} package was also briefly compared to our approach. It seems that the \texttt{pseudo} function provides very similar results but with a faster computational time. However, since there is no available information on the implementation of the \texttt{pseudo} function, it is not possible to clearly evaluate which of the two formulas (the one obtained with the infinitesimal jackknife and the one obtained from the Von Mises formula) has the smallest computational complexity.

\begin{table}[ht]
\centering
\begin{tabular}{|c|cccc|cccc|}
  \hline
  & \multicolumn{4}{c|}{Jackknife}& \multicolumn{4}{c|}{Approximated formula}\\
 $n$ & Bias($\hat\beta$) & SE($\hat\beta$) & MSE($\hat\beta$) & Time & Bias($\hat\beta$) & SE($\hat\beta$) & MSE($\hat\beta$) & Time \\ \hline
100 & -0.006 & 0.273 & 0.075 & 0.211 s& -0.004 & 0.269 & 0.073 & 0.015 s\\ 
  & -0.020 & 0.368 & 0.136 &  & -0.022 & 0.361 & 0.131 &  \\ 
   & -0.003 & 0.368 & 0.135 &  & -0.006 & 0.361 & 0.130 &  \\ 
   & -0.009 & 0.364 & 0.132 &  & -0.013 & 0.357 & 0.128 &  \\ 
 500 & 0.007 & 0.115 & 0.013 & 1.490 s& 0.007 & 0.115 & 0.013 & 0.054 s\\ 
   & -0.005 & 0.152 & 0.023 &  & -0.005 & 0.151 & 0.023 &  \\ 
   & -0.010 & 0.157 & 0.025 &  & -0.010 & 0.156 & 0.024 &  \\ 
   & -0.003 & 0.154 & 0.024 &  & -0.004 & 0.153 & 0.024 &  \\ 
  1,000 & 0.002 & 0.081 & 0.007 & 4.084 s& 0.002 & 0.081 & 0.007 & 0.163 s\\ 
   & -0.003 & 0.110 & 0.012 &  & -0.003 & 0.110 & 0.012 &  \\ 
   & -0.004 & 0.113 & 0.013 &  & -0.004 & 0.112 & 0.013 &  \\ 
   & -0.003 & 0.107 & 0.011 &  & -0.003 & 0.107 & 0.011 &  \\ 
 10,000 & 0.002 & 0.026 & 0.001 & 4.429 min & 0.002 & 0.026 & 0.001 & 14.194 s\\ 
  & -0.001 & 0.037 & 0.001 &  & -0.001 & 0.037 & 0.001 &  \\ 
   & 0.001 & 0.036 & 0.001 &  & 0.001 & 0.036 & 0.001 &  \\ 
   & -0.001 & 0.034 & 0.001 &  & -0.001 & 0.034 & 0.001 &  \\ 
   \hline
\end{tabular}
\caption{Simulation results for the estimation
       of $\beta$ in the RMST model~\eqref{simuRC} based on pseudo-regression with the Kaplan-Meier estimator on $33\%$ of right-censored data. In the pseudo-regression, 
       the true jackknife is compared to the approximated pseudo-estimates.} \label{simu:RC}
\end{table}

\subsection{Interval-censored data}\label{sec:simuIC}

For interval censored data the survival function is estimated from the pch model, as detailed in Section~\ref{sec:param_model}. Using this model, estimation of the model parameter $\alpha_0$ is performed using the EM algorithm, as presented in~\cite{bouaziz21}. An alternative method could be to directly maximise the observed likelihood but this would result in implementing the Newton-Raphson algorithm for each jackknife sample with inversion of a Hessian matrix of full rank which, in turn, would result in unstable results. In the EM algorithm, the M-step is explicit and as a result the computation of the jackknife methods is always stable. We refer the reader to~\cite{bouaziz21} for more details on the two methods. The approximated method is implemented from the result in Proposition~\ref{paramsurv} and details on the computation of the score vector and Hessian matrix are detailed in Section~\ref{sec:scoreHessian_pch} of the Appendix.



We assume Model~\eqref{simuRC} with the same values of $\sigma$ and $\tau$. Then, in order to simulate interval-censored data, a total of $K=5$ visits were simulated such that $V_1\sim\mathcal U[0,6]$ and $V_k=V_{k-1}+ U[0,2]$, for $k=2,\ldots,K$. The observations for which $T^*_i<V_1$ correspond to left-censored observations with $L_i=0$ and $R_i=V_1$, the observations for which $T^*_i>V_K$ correspond to right-censored observations with $L_i=V_K$ and $R_i=\infty$, and the observations for which $V_{k-1}<T^*_i<V_{k}$ ($k=2,\ldots, K$) correspond to strictly interval-censored observations with $L_i=V_{k-1}$ and $R_i=V_{k}$. This resulted in $14.6\%$ of left-censored data, $52.07\%$ of interval-censored data and $33.33\%$ of right-censored data. For interval-censored data, the average length of the intervals was approximately equal to $1.34$. The pch model with cuts equal to $4, 5, 6, 7$ was used for the computation of the survival estimator. 
 The pseudo-observations were generated based on the standard jackknife and on the approximated formulas and the results for the RMST model are presented in Table~\ref{simu:IC}. 

Again the results between the jackknife and the approximate formula are almost identical while there is a huge gain in terms of computational time for the approximated formula. The approximated formula is $107$, $198$ and $310$ times faster than the jackknife method for $n=200$, $n=500$ and $n=1,000$ respectively. It should be noted that the cuts must be carefully chosen in the pch model. In particular, the regularity conditions of Equation~\eqref{eq:condIC} must be satisfied. If there are only few values of $L_i$ and $R_i$ that intersect a cut $[c_{k-1},c_k]$ or if the proportion of $L_i$'s such that $L_i>c_{k-1}$ is too low then the pseudo-values can be incorrect (both for the jackknife method or using our approximated formula) which will in turn result in a poor performance of the parameters estimation. On the other hand, if the regularity conditions hold, the choice of the cuts will only have a minor impact on the performance of the estimator of $\beta_0$ and will lead to similar results. A supplementary simulation scenario for interval-censored data with $\tau$ equal to infinity is also presented in the Supporting Information. 

\begin{table}[ht]
\centering
\begin{tabular}{|c|cccc|cccc|}
  \hline
  & \multicolumn{4}{c|}{Jackknife}& \multicolumn{4}{c|}{Approximated formula}\\
 $n$ & Bias($\hat\beta$) & SE($\hat\beta$) & MSE($\hat\beta$) & Time & Bias($\hat\beta$) & SE($\hat\beta$) & MSE($\hat\beta$) & Time \\ \hline
200 & -0.169 & 0.220 & 0.077 & 42.047s & -0.167 & 0.219 & 0.076 & 0.392s \\ 
   & 0.019 & 0.310 & 0.096 &  & 0.018 & 0.308 & 0.095 &  \\ 
  & 0.015 & 0.305 & 0.094 &  & 0.013 & 0.304 & 0.093 &  \\ 
   & 0.057 & 0.293 & 0.089 &  & 0.055 & 0.293 & 0.089 &  \\ 
  500 & -0.181 & 0.141 & 0.053 & 2.408 min & -0.181 & 0.140 & 0.052 & 0.731s \\ 
   & 0.036 & 0.190 & 0.037 &  & 0.036 & 0.190 & 0.037 &  \\ 
   & 0.035 & 0.191 & 0.038 &  & 0.035 & 0.191 & 0.038 &  \\ 
   & 0.077 & 0.182 & 0.039 &  & 0.076 & 0.182 & 0.039 &  \\ 
  1,000 & -0.184 & 0.101 & 0.044 & 6.675 min & -0.183 & 0.101 & 0.044 & 1.292s \\ 
   & 0.035 & 0.136 & 0.020 &  & 0.035 & 0.136 & 0.020 &  \\ 
   & 0.033 & 0.142 & 0.021 &  & 0.032 & 0.142 & 0.021 &  \\ 
   & 0.072 & 0.130 & 0.022 &  & 0.072 & 0.130 & 0.022 &  \\ 
   \hline
\end{tabular}
\caption{Simulation results for the estimation of 
$\beta$ in the RMST model~\eqref{simuRC} based on pseudo-regression with $14.6\%$ of left-censored data, $52.07\%$ of 
interval-censored data and $33.33\%$ of right-censored data. The piecewise constant hasard model with cuts equal 
to $4, 5, 6, 7$ was used for the estimation of the survival function in the computation of the pseudo-observations. 
In the pseudo-regression, the true jackknife is compared to the approximated pseudo-estimates.} \label{simu:IC}
\end{table}

\section{Illustrative real data examples}\label{sec:realdata}

\subsection{The Cardiovascular Health Study (CHS)}
 
In this data example, we mimic the analysis of the Cardiovascular Health Study (CHS) as it was performed in~\cite{zhao2021deep}. This study was initiated in $1987$ to determine the risk factors for development and progression of cardiovascular disease (CVD) in older adults. The event of interest was time to CVD. In~\cite{zhao2021deep}, the author considers a subsample of $5,380$ individuals of whom $65.2\%$ had CVD during the study period and the others were right-censored. The aim of the study was to estimate the conditional RMST with $29$ covariates and $\tau=5$ years.

The methodology proposed in~\cite{zhao2021deep} uses pseudo-observations and implements a deep neural network directly on the pseudo-observations of the RMST, that is the $g$ link function presented in Section~\ref{sec:backgrounds} is a neural network. Moreover, a training dataset including $75\%$ of the observations and a test set based on the remaining $25\%$ of the data are built in order to evaluate the prediction performance of the method. This split of the data between training and test sets is repeated $10$ times. At each repetition, the pseudo-observations must be entirely computed but only on the training datasets. This results in computing the pseudo-observations for the RMST for $10$ samples of size $4,035$. We computed those pseudo observations from the jackknife method and the approximated formula. The former was computed in approximately $21.3$ seconds while the latter took $1.9$ seconds. Therefore, our approximated formula is more than $11$ times faster than the original jackknife method. Since building a neural network is computationally expensive and needs to be implemented for all the training samples, this reduction in computation time is a major advantage for our approximated formula. Of note, the results of the analysis implemented with the approximated formula are identical to the original analysis (based on the jackknife method) and are therefore omitted.

\subsection{The Signal $\text{Tandmobiel}^{\Circled{\text{R}}}$ data}

In this section, our aim is to analyse the Signal $\text{Tandmobiel}^{\Circled{\text{R}}}$ data using three different models: the standard Cox model, a logistic model for the conditional survival function and the conditional RMST model presented in Equation~\eqref{eq:RMST}. This dataset is part of the \texttt{bayesSurv} R package. Those data were collected from a longitudinal dental survey of $4,468$ 
school children born in $1989$, who were annually examined by a dentist. The time scale is age in years. The dataset is composed of $0.68\%$ of left-censored data, $61.69\%$ of strictly interval-censored data and $37.63\%$ of right-censored data. Our aim is to study the emergence of the tooth number $14$ which is a permanent first premolar. The covariates used for the analysis are: gender (binary variable equal to $1$ for boys, $0$ for girls) and the number of decayed or missing deciduous first molars due to caries among teeth $54$, $64$, $74$, $84$ of the dataset. This covariate is thus discrete taking values between $0$ and $4$. These data were previously studied by~\cite{lesaffre2005overview} using the Accelerated Failure Time model (AFT). In our analysis, the survival function is estimated from the pch model using the whole dataset and the pseudo-observations are then computed from the approximated formula in Proposition~\ref{paramsurv}. There are $126$ individuals with missing covariates and the generalised estimating equation used to implement our models is therefore applied to this reduced dataset composed of $4,342$ pupils. 

\begin{figure}[!htb]
\begin{tabular}{cc}
\includegraphics[width=0.47\textwidth,height=0.4\textwidth]{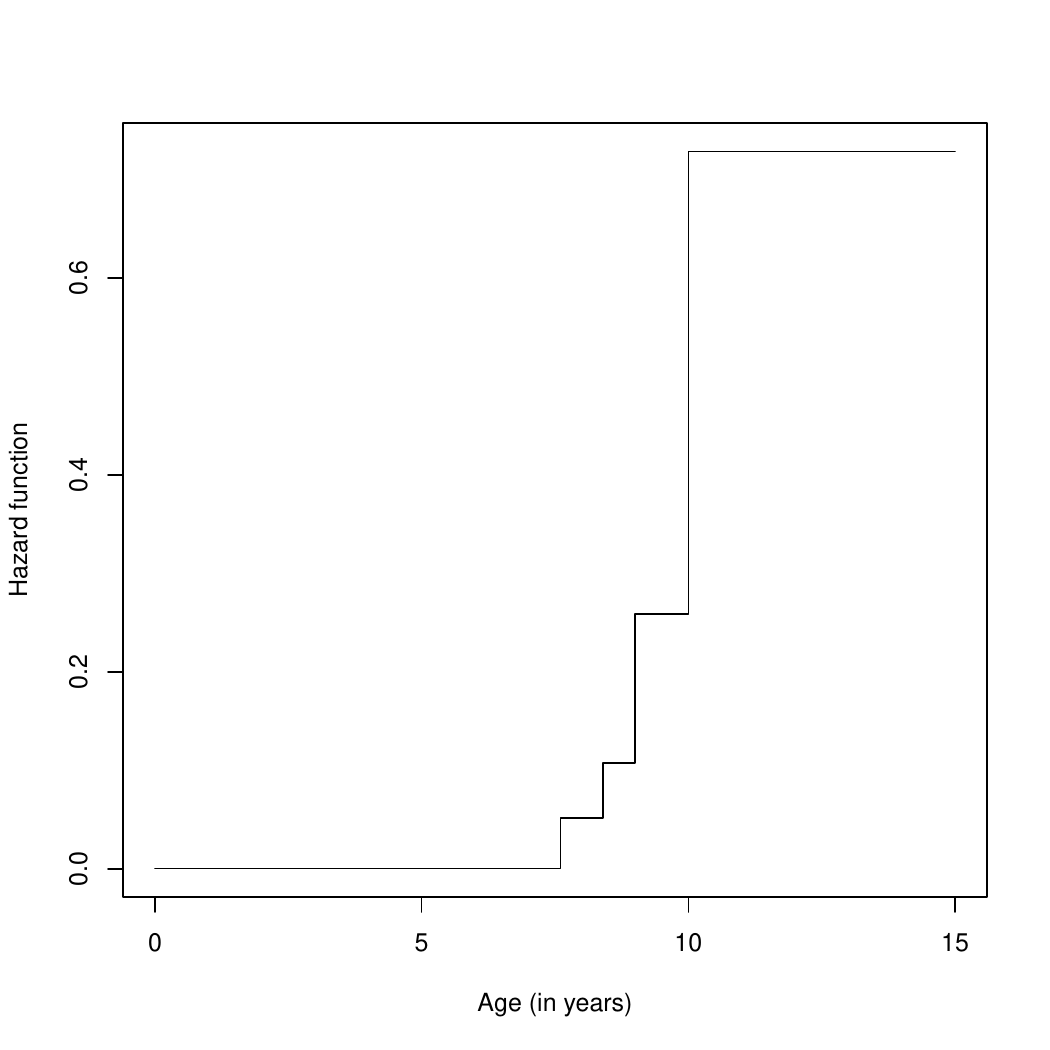}&\includegraphics[width=0.47\textwidth,height=0.4\textwidth]{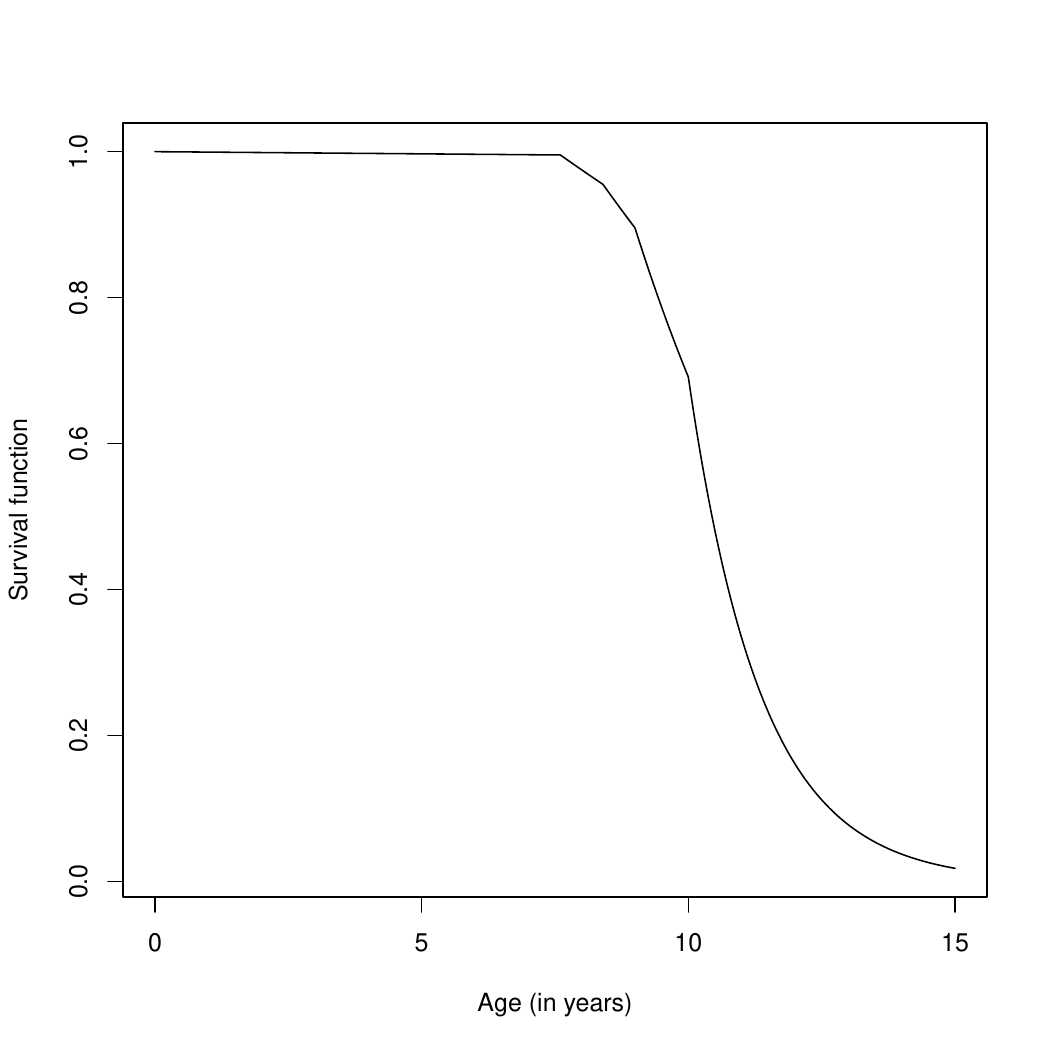}
\end{tabular}
\caption{Distribution of time to emergence of the tooth number $14$. On the left: estimated hazard function. On the right: estimated survival function. Those estimates were obtained from the pch model with cuts equal to $7.6, 8.4, 9$ and $10$.}
\label{fig:haz_surv_dental}
\end{figure}

In the pch model, the number of cuts and locations were chosen using the adaptive-ridge algorithm developed in~\cite{bouaziz21}. This led to the selection of the four cuts $7.6$, $8.4$, $9$ and $10$. Since the maximum likelihood estimator has converged this entails that the regularity conditions of Equation~\eqref{eq:condIC} are satisfied. We can also easily check them empirically: in particular there are $3\%$ of strictly interval censored observations whose left intervals fell before $7.6$, $40\%$ of left intervals that fell after $10$ and the percentage of strictly interval censored observations that intersect each other is high (values not shown). The corresponding estimated hazard and survival functions are displayed in Figure~\ref{fig:haz_surv_dental}. We observe a low estimated hazard value (equal to $6\!\cdot \! 10^{-4}$) from age $0$ until age $7.6$ due the low percentage of left intervals that fell before $7.6$. This yields a very flat decay of the survival function on this time period, then the decay increases drastically for the four other time periods $[7.6,8.4]$, $[8.4,9]$, $[9,10]$ and $[10,\infty)$. For illustration, we estimate from the survival function that approximately $83.39\%$ of the teeth will emerge between age $7.6$ and $12$. 

In order to implement a Cox model from the pseudo-observations, we need to set a grid of time points on which the baseline hazard rate is estimated (see~\cite{andersen2003generalised}). We choose the time points $t_1=8, t_2=9, t_3=10, t_4=11, t_5=12$ and we use the link function $g(x)=\log(-\log(x))$ in the pseudo-regression which leads to the Cox model:
\begin{align*}
\log(-\log(S(t_m\mid Z_i)) =\log \Lambda_0(t_m) + Z_i^{\top}\beta, \quad m=1,\ldots, 5,
\end{align*}
where $\Lambda_0(t_m)$ represents the cumulative baseline hazard function at time $t_m$. The results are presented in Table~\ref{tab:Cox_teeth}. The second column provides the cumulative baseline hazard at different time points and the hazard ratios of the two covariates. The last column displays the Wald tests which are all extremely significant. We clearly see that the hazard for the emergence of the tooth is an increasing function of time. Also, under the proportional hazards assumption, boys have an increased risk as compared to girls with an hazard ratio equal to $1.47$ and the hazard ratio for one supplementary decayed or missing deciduous first molar due to caries equals $1.13$.

For the logistic regression model, we study the probability for the emergence of the tooth before a fixed time point. We conduct two separate analysis, one for the time point $9$ and another one for the time point $12$ which correspond to the $11\%$ and $84\%$ estimated quantiles of $T^*$, respectively, thus corresponding to early and late emergence of the tooth. There is no guarantees that the pseudo-observations are in the interval $[0,1]$ and we therefore set the negative values to $0$ and the values greater than $1$ to $1$. The model is the following:
\begin{align*}
\text{logit}(1-S(t\mid Z_i))=\gamma+Z_i^{\top}\beta,
\end{align*}
where $t$ is either equal to $9$ or $12$, $\gamma$ is the intercept and $\beta$ is a two dimensional vector representing the effect of the covariates and logit is the classical logistic function ($\text{logit}(x)=\log(x/(1-x))$. It should be noted that the model can be directly implemented with the \texttt{glm} or \texttt{geese} functions by using one minus the pseudo-observations as input. The results are presented in Table~\ref{tab:Logistic_teeth}. We first observe that the effect of the two covariates are highly significant but they differ depending on the time endpoint: the odds ratio for the number of decayed or missing deciduous first molars decreases from $\exp(0.2808)\approx 1.32$ before age $9$ to $\exp(0.1080)\approx 1.11$ before age 12 while the odds ratio for the gender effect increases from $\exp(0.2978)\approx 1.35$ before age $9$ to $\exp(0.5284)\approx 1.70$ before age $12$. This strongly suggests that the number of decayed or missing deciduous first molars is mostly responsible for the early teeth emergence while gender (with boys having a higher risk) is mostly responsible for the late tooth emergence. It is also interesting to compare the probabilities of the tooth emergence before age $9$ between girls with no decayed or missing deciduous first molars ($5.49\%$), boys with no decayed or missing deciduous first molars first molars ($7.26\%$) and boys with $4$ decayed or missing deciduous first molars ($19.39\%$). We can similarly compare the probabilities of the tooth emergence before age $12$ between girls with no decayed or missing deciduous first molars ($70.63\%$), boys with no decayed or missing deciduous first molars ($80.31\%$) and boys with $4$ decayed or missing deciduous first molars ($86.27\%$).


Finally, two RMST analysis were conducted with $\tau=9$ and $\tau=12$.
The estimated regression parameters in the RMST model along with their Wald test are presented in Table~\ref{tab:RMST_teeth}. For $\tau=9$ we observe a weak effect of the covariates with an intercept that is almost equal to $\tau$, highlighting that most emergences of the tooth will occur after $9$ years of age. As a matter of fact, gender is not significant and the number of decayed or missing deciduous first molars is highly significant but with a weak effect. The number of decayed or missing deciduous first molars will accelerate the emergence of the tooth with $1$ decayed molar (respectively $4$ decayed molars) yielding a reduction of $0.0097$ years (respectively $0.0390$ years) for the emergence of the tooth. For $\tau=12$ the effect of gender is now highly significant, meaning that gender only plays a role for late emergence of the tooth (a finding that was also observed with the logistic regression model). The emergence of the tooth for boys arrives on average $0.3336$ years earlier than for girls. The number of decayed or missing deciduous first molars is also highly significant with $1$ decayed molar (respectively $4$ decayed molars) yielding a reduction of $0.1303$ years (respectively $0.5211$ years) for the emergence of the tooth. 

We also tried to repeat the procedure using different cut values in the pch model and as already observed in the simulation study, this led to very similar results, for all three models. The results from~\cite{lesaffre2005overview} obtained using the AFT models were similar to our findings except that the authors did not provide statistical tests for the effects of the covariates and it was not possible from their method to detect that gender had mainly a role for late emergence of the tooth.

Finally, based on the approximated formulas developed in this paper, the whole procedure (computation of the pseudo-observations and implementation of the generalised estimating equations) took about $1.78$ seconds for the RMST analysis (the computation times are similar for the other two models). The method was not implemented using the classical jackknife method but according to the simulation study it would have taken $29.8$ minutes to obtain the pseudo-observations, since in the simulation study the time for the jackknife procedure was evaluated at 6.7 minutes for $n=1,000$ (see Sections~\ref{sec:simuIC}). Also, the results would have been identical, thus highlighting the relevance of the proposed approach in practical situations.



\begin{table}[htb!]
	\centering
	\begin{tabular}{|c|c|c|c|c| }
		\hline
		\textbf{Covariates}&\textbf{effect}&\textbf{exp. effect}& \textbf{se}&\textbf{p-value}\\\hline
  Intercept at time $8$ & -4.4510 & 0.0117 & 0.1175 & $<10^{-15} $\\ 
  Intercept at time $9$ & -2.6999 & 0.0672 & 0.0727 & $<10^{-15}$ \\ 
  Intercept at time $10$ & -1.4461 & 0.2355 & 0.0463 & $<10^{-15} $\\ 
  Intercept at time $11$ & -0.3206 & 0.7257 & 0.0358 & $<10^{-15}$ \\ 
  Intercept at time $12$ & 0.2293 & 1.2577 & 0.0340 & $<10^{-10}$\\ 
  Gender ($1=$ boy) & 0.3885 & 1.4747 & 0.0395 & $<10^{-15}$\\ 
  Nb of decayed molars & 0.1249 & 1.1330 & 0.0130 & $<10^{-15} $\\ 
		\hline
	\end{tabular}
	\caption{Cox Model for the time to emergence of the tooth $14$ with the covariates gender and number of decayed or missing deciduous first molars due to caries among teeth $54$, $64$, $74$, $84$. The time has been discretised at $t=8, 9, 10,11,12$ in the pseudo-regression approach. \textbf{se} represents the standard estimate of the regression parameter.}\label{tab:Cox_teeth}
\end{table}

\begin{table}[htb!]
	\centering
	\begin{tabular}{|c|c|c|c|c|c|c| }
	 \hline
  &\multicolumn{3}{c| }{$\text{time}=9$}& \multicolumn{3}{c| }{$\text{time}=12$}\\ 
		\hline
		\textbf{Covariates}& \textbf{effect} & \textbf{se} & \textbf{p-value} & \textbf{effect} & \textbf{se} & \textbf{p-value} \\ 
  \hline
Intercept & -2.8458 & 0.0761 & $<10^{-15} $ & 0.8777 & 0.0501 & $<10^{-15} $ \\ 
  Gender ($1=$ boy) & 0.2978 & 0.0819 & 0.0003 & 0.5284 & 0.0599 & $<10^{-15} $ \\ 
  Nb of decayed molars & 0.2808 & 0.0260 & $<10^{-15} $ & 0.1080 & 0.0194 & $2.7 \times 10^{-8}$ \\ 
		\hline
	\end{tabular}
	\caption{Two logistic models for the probability that the emergence of the tooth $14$ occurred before time $9$ and $12$ with respect to the covariates gender and number of decayed or missing deciduous first molars due to caries among teeth $54$, $64$, $74$, $84$. \textbf{se} represents the standard estimate of the regression parameter.}\label{tab:Logistic_teeth}
\end{table}

\begin{table}[htb!]
	\centering
	\begin{tabular}{|c|c|c|c|c|c|c| }
	 \hline
  &\multicolumn{3}{c| }{$\tau=9$}& \multicolumn{3}{c| }{$\tau=12$}\\ 
		\hline
		\textbf{Covariates}&\textbf{effect}& \textbf{se}&\textbf{p-value}&\textbf{effect}& \textbf{se}&\textbf{p-value}\\\hline
		Intercept &$8.9851$ & $0.0047$ & $<10^{-15}$&$10.8755$ & $0.0306$ & $<10^{-15}$\\
		Gender ($1=$ boy) &$-0.0097$ & $0.0066$ & $0.1422$&$-0.3336$ & $0.0361$ & $<10^{-15}$\\
		Nb of decayed molars & $-0.0180$ & $0.0024$ & $8.1379\times 10^{-14}$& $-0.1303$ & $0.0120$ & $< 10^{-15}$\\
		\hline
	\end{tabular}
	\caption{Restricted Mean Survival Time Model for the time to emergence of the tooth $14$ with the covariates gender and number of decayed or missing deciduous first molars due to caries among teeth $54$, $64$, $74$, $84$. Two values of $\tau$ are analysed in Equation~\eqref{eq:RMST}. \textbf{se} represents the standard estimate of the regression parameter.}\label{tab:RMST_teeth}
\end{table}


\section{Conclusion}

In this paper, we presented asymptotic formulas for computing pseudo-observations for time to event data. In the context of right-censored data, those formulas are based on the Kaplan-Meier estimator of the survival function. When dealing with interval-censoring our formulas were developed for a general class of parametric models. Pseudo-regression is an appealing tool when the goal is to directly model a complex quantity of interest, such as the RMST, cumulative incidence functions in a competing risk setting or transition probabilities for multi-state models. Our formulas were precisely developed for the RMST but they could be easily extended for those other quantities of interest. 
While the pseudo-values approach is originally based on the jackknife procedure, our formulas only involve quantities computed on the initial sample. 
This results in a drastic reduction of the computational time, which is an interesting feature when dealing with large dataset or when the data are interval-censored, since in that case, the estimators are computationally intensive.

There has been an increasing interest of the pseudo-values approach in the machine learning community. After having computed the pseudo-observations, standard machine learning models can be applied to those new observations, by simply ignoring the censoring. In particular, this methodology has been applied for estimating the survival function in~\cite{zhao2019dnnsurv},~\cite{zhao2020deep},~\cite{feng2021bdnnsurv} or the RMST in~\cite{zhao2021deep} based on neural networks that were directly applied on the pseudo-observations. Therefore, our formulas are particularly interesting in those settings where the dataset can be extremely large and the algorithm usually relies on a cross-validation procedure. Using our approximated formulas results in a significant gain in terms of computation time as illustrated on the real data analysis. Also, the approximations made by our formulas are extremely precise, even for moderate sample sizes, as shown in the simulation study. Surprisingly, we also saw that our formulas are more robust than the original jackknife method which sometimes fails due to some rare extreme values. For all these reasons, we advocate the use of our asymptotic formulas in practical situations.

As a reviewer pointed out, there exists an alternative to our fast formula in the case of right-censored data, in the \texttt{pseudo} function of the \texttt{survival} package. It is based on the infinitesimal jackknife and can be applied to the Kaplan-Meier estimator or to the Breslow estimator (defined as the exponential of minus the Nelson-Aalen estimator). A theoretical comparison of those approximations with the Von-Mises approach developed in this paper has been made in the Supporting Information. It is shown that all three approximations are very similar for small sample sizes and asymptotically equivalent. It would now be interesting to extend the infinitesimal jackknife to the parametric setting, in order to apply it to interval-censored data, as there already exists some approaches based on the infinitesimal jackknife to compute delta-beta residuals, which are very closely related to pseudo-observations. 

Finally, new theoretical results for parametric pseudo-regression were also developed in this paper. Those results indicate that parametric pseudo-regression is only valid up to two extra terms but simulation studies suggest that they are reasonably small in practice and do not significantly impact the performance of the final estimator, as long as the number of covariates is not too large. They also suggest to use flexible parametric models such as the pch model. As an extension it would be interesting to study the theoretical validity of splines models for pseudo-regression such as the method developed in~\cite{nygaard2020regression}. This is let to future research.

\section{Software}

The asymptotic formulas developed in this paper for the pseudo-values of the survival function and the RMST can be implemented using the GitHub package \texttt{FastPseudo} available at \href{https://github.com/obouaziz/FastPseudo}{https://github.com/obouaziz/FastPseudo}. The package can deal with both right-censored or interval-censored data. In the latter case, the formulas are implemented for the pch model.





\section*{Acknowledgments}

We thank the reviewers for their very constructive criticisms and comments that have helped improve the paper. We also thank Per Kragh Andersen and Terry Therneau for their very interesting discussions about the connection between the infinitesimal jackknife and the Von Mises approximation for right-censored data.

\section{Appendix}

\subsection{Proof of Proposition~\ref{KM}}

Introduce $H_1(t)=\mathbb P(T\leq t,\Delta=1)$, $H_0(t)=\mathbb P(T\leq t,\Delta=0)$, $H(t)=\mathbb P(T\geq t)$ and their empirical counterparts, $\hat H_1(t)=\sum_i I(T_i\leq t,\Delta_i=1)/n$, $\hat H_0(t)=\sum_i I(T_i\leq t,\Delta_i=0)/n$, $\hat H(t)=\sum_i I(T_i\geq t)/n$. 
Let $\psi(A)(s,t]=\Prodi_{s<u\leq t}(1+dA(u))$, such that $\psi(\Lambda)(0,t]=S(t)$, where $\Lambda(t)$ is the cumulative hazard function and $\psi(\hat \Lambda)(0,t]=\hat S(t)$. We have the following Von-Mises expansion \citep*[see][]{gill1994lectures,van1996}:
\begin{align*}
\hat S^{(-l)}(t) &= \hat S(t)-\hat S(t)(\hat \Lambda^{(-l)}(t)-\hat \Lambda(t))+o_{\mathbb P}(\hat \Lambda^{(-l)}(t)-\hat \Lambda(t)).
\end{align*}
We now derive a Von-Mises expansion for $\hat \Lambda^{(-l)}(t)-\hat \Lambda(t)$. The cumulative hazard function and its estimator can be defined as functions of $H$, $H_1$ and of $\hat H_1$, $\hat H$ respectively where $\Lambda(t)=g(H_1,H):=\int_0^t dH_1(u)/H(u)$ and $\hat \Lambda(t)=g(\hat H_1,\hat H)$. 
We have the following Von-Mises expansion:
\begin{align*}
\hat \Lambda^{(-l)}(t) &= \hat \Lambda(t) + g'_{(\hat H_1,\hat H)}(\hat H_1^{(-l)}-\hat H_1,\hat H^{(-l)}-\hat H)+o_{\mathbb P}(n^{-1}),
\end{align*}
where $g'$ is the Hadamard derivative of $g$, which is equal to \citep[see][]{gill1994lectures, van1996}:
\begin{align*}
g'_{(H_1,H)}(h_1,h) & = \int_0^t \frac{dh_1}{H}- \int_0^t \frac{h_2dH_1}{H^2}\cdot
\end{align*}
The $o_{\mathbb P}(n^{-1})$ term above comes from the expressions:
\begin{align*}
\hat H_1^{(-l)}(t)-\hat H_1(t) & = \frac 1{n(n-1)} \sum_{i=1}^n I(T_i\leq t,\Delta_i=1)-\frac{I(T_l\leq t,\Delta_l=1)}{n-1}
\end{align*}
and
\begin{align*}
\hat H^{(-l)}(t)-\hat H(t) & = \frac 1{n(n-1)} \sum_{i=1}^n I(T_i\geq t)-\frac{I(T_l\geq t)}{n-1},
\end{align*}
which entail as a consequence that $\hat H_1^{(-l)}(t)-\hat H_1(t)$ and $\hat H^{(-l)}(t)-\hat H(t)$ are $O_{\mathbb P}(n^{-1})$.
Moreover, using those expressions we have
\begin{align*}
g'_{(\hat H_1,\hat H)}(\hat H_1^{(-l)}-\hat H_1,\hat H^{(-l)}-\hat H) & = \frac 1{n-1}\int _0^{t}\frac{d\hat H_1(u)}{\hat H(u)}-\frac 1{n-1}\frac{I(T_l\leq t,\Delta_l =1)}{\hat H(T_l)}\\
&\quad-\frac 1{n-1}\int _0^{t}\frac{d\hat H_1(u)}{\hat H(u)}+\frac 1{n-1}\int_0^t\frac{I(T_l\geq u)d\hat H_1(u)}{(\hat H(u))^2},\\
&=-\frac 1{n-1}\frac{I(T_l\leq t,\Delta_l =1)}{\hat H(T_l)}+\frac 1{n-1}\int_0^t\frac{I(T_l\geq u)d\hat H_1(u)}{(\hat H(u))^2}\cdot
\end{align*}
Gathering all the different parts, we obtain
\begin{align*}
\hat S_{(l)}(t)&=\hat S(t)+\hat S(t)\left( \int _0^{T_l\wedge t}\frac{d\hat H_1(u)}{(\hat H(u))^2}-\frac{I(T_l\leq t,\Delta_l =1)}{\hat H(T_l)}\right)+o_{\mathbb P}(1)\\
&=\hat S(t)+\hat S(t)\left( \int _0^{t}\frac{I(u\leq T_l)}{\hat H(u)}d\hat \Lambda(u)-\int_0^t\frac{dN_l(u)}{\hat H(u)}\right)+o_{\mathbb P}(1)\\
&=\hat S(t)-\int_0^t \frac{d\hat M_l(u)}{\hat H(u)}+o_{\mathbb P}(1)\cdot
\end{align*}
The approximation for the RMST is then obtained by directly integrating the previous equation as it actually follows that the convergence holds in the Skorohod space $D[0,\tau]$~\citep*[see][]{gill1994lectures,van1996} and therefore the convergence holds uniformly with respect to $t\in[0,\tau]$. 

\subsection{Proof of Proposition~\ref{paramsurv}}

Starting with Equation~\eqref{eq:pseudo_param} we will first prove that $\sqrt n\, \varepsilon_n-\sqrt{n-1}\, \varepsilon^{(-l)}_n$ tends to $0$ in probability as $n$ tends to infinity. Set
\begin{align*}
\tilde I_n & =-\frac 1n \sum_{i=1}^n \nabla^2 \log f(X_i;\tilde{\alpha}), \quad \tilde I_n^{(-l)} =-\frac 1{n-1} \sum_{i\neq l} \nabla^2 \log f(X_i;\tilde{\alpha}), 
\end{align*}
where $\tilde{\alpha}$ lies between $\hat\alpha$ and $\alpha_0$. We have:
\begin{align*}
&\sqrt n\, \varepsilon_n-\sqrt{n-1}\, \varepsilon^{(-l)}_n\\
&\quad=\left(\tilde I_n^{-1}-I^{-1}\right)\sum_{i=1}^n \nabla \log f(X_i;\alpha_0)- \left(\Big(\tilde I_n^{(-l)}\Big)^{-1}-I^{-1}\right)\sum_{i\neq l}^n \nabla \log f(X_i;\alpha_0)\\
&\quad = \left(\Big(\tilde I_n^{(-l)}\Big)^{-1}-I^{-1}\right)\nabla \log f(X_l;\alpha_0)+ \left(\tilde I_n^{-1}-\Big(\tilde I_n^{(-l)}\Big)^{-1}\right)\sum_{i=1}^n \nabla \log f(X_i;\alpha_0).
\end{align*}
Clearly for $I$ positive definite, $ \left(\Big(\tilde I_n^{(-l)}\Big)^{-1}-I^{-1}\right)\nabla \log f(X_l;\alpha_0)$ tends to $0$ in probability. Since $\sum_{i=1}^n  \nabla \log f(X_i;\alpha_0)/n$
tends to $\mathbb E(\nabla \log f(X;\alpha_0))=0$ in probability, we just need to prove that $\tilde I_n^{-1}-\Big(\tilde I_n^{(-l)}\Big)^{-1}=O_{\mathbb P}(1/n)$ to conclude the proof. Write:
\begin{align*}
\tilde I_n^{-1}-\Big(\tilde I_n^{(-l)}\Big)^{-1} & =\tilde I_n^{-1}(\tilde I_n^{(-l)}-\tilde I_n)\Big(\tilde I_n^{(-l)}\Big)^{-1}.
\end{align*}
From the law of large numbers, $\tilde I_n^{-1}$ and $\Big(\tilde I_n^{(-l)}\Big)^{-1}$ tend towards $I^{-1}$ in probability and
\begin{align*}
\tilde I_n^{(-l)}-\tilde I_n & =-\frac 1{n(n-1)}\sum_{i\neq l} \nabla^2 \log f(X_i;\tilde \alpha)+\frac 1n \nabla^2 \log f(X_l;\tilde \alpha)= O_{\mathbb P}\left(1/n\right).
\end{align*}
This proves that
\begin{align}\label{eq:alpha}
n\hat\alpha-(n-1) \hat\alpha^{(-l)} = \alpha_0+I^{-1}\nabla \log f(X_l;\alpha_0)+o_{\mathbb P}(1).
\end{align}
Using the consistency of $\Lambda(t;\hat\alpha)$ towards $\Lambda(t;\alpha_0)$ from standard maximum likelihood theory, we now write a Taylor expansion for the cumulative hazard function around $\alpha_0$:
\begin{align}\label{eq:Lambda}
\Lambda(t;\hat\alpha)=\Lambda(t;\alpha_0)+(\hat\alpha-\alpha_0)^{\top}\nabla\Lambda(t;\alpha_0)+\frac 12 (\hat\alpha-\alpha_0)^{\top}\nabla^2\Lambda(t;\tilde\alpha)(\hat\alpha-\alpha_0),
\end{align}
where $\tilde\alpha$ lies between $\hat\alpha$ and $\alpha_0$. We also write a Taylor expansion for the function $x\mapsto \exp(-x)$ around $0$:
\begin{align*}
\exp(-(\Lambda(t;\hat\alpha)-\Lambda(t;\alpha_0))) & =1-\big(\Lambda(t;\hat\alpha)-\Lambda(t;\alpha_0)\big)+\frac 12 \big(\Lambda(t;\hat\alpha)-\Lambda(t;\alpha_0)\big)^2\\
&\quad-\frac 16 e^{\xi_n}\big(\Lambda(t;\hat\alpha)-\Lambda(t;\alpha_0)\big)^3,
\end{align*}
with $\xi_n$ tends to $0$ in probability as $n$ tends to infinity. This can be rewritten as:
\begin{align*}
S(t;\hat\alpha) & = S(t;\alpha_0)+S(t;\alpha_0)\left(-\big(\Lambda(t;\hat\alpha)-\Lambda(t;\alpha_0)\big)+\frac 12 \big(\Lambda(t;\hat\alpha)-\Lambda(t;\alpha_0)\big)^2\right)+o_{\mathbb P}(1/n),
\end{align*} 
using the fact that $\sqrt n (\Lambda(t;\hat\alpha)-\Lambda(t;\alpha_0))$ converges in distribution towards a centred gaussian variable with finite variance from standard results on maximum likelihood theory and the delta-method.
As a result,
\begin{align}\label{eq:S}
n S(t;\hat\alpha) - (n-1) S(t;\hat\alpha^{(-l)}) &=S(t;\alpha_0)+A_{n,1}+A_{n,2}+o_{\mathbb P}(1),
\end{align}
where
\begin{align*}
A_{n,1}&=-S(t;\alpha_0)\left(n\big(\Lambda(t;\hat\alpha)-\Lambda(t;\alpha_0)\big)-(n-1)\big(\Lambda(t;\hat\alpha^{(-l)})-\Lambda(t;\alpha_0)\big)\right),\\
A_{n,2}&=S(t;\alpha_0) \left(n\big(\Lambda(t;\hat\alpha)-\Lambda(t;\alpha_0)\big)^2-(n-1)\big(\Lambda(t;\hat\alpha^{(-l)})-\Lambda(t;\alpha_0)\big)^2\right)\frac {1}2\cdot
\end{align*}
We start with the $A_{n,2}$ term. 
From Equation~\eqref{eq:Lambda} we have:
\begin{align*}
n\big(\Lambda(t;\hat\alpha)-\Lambda(t;\alpha_0)\big)^2&= n(\hat\alpha-\alpha_0)^{\top}\nabla\Lambda(t;\alpha_0)\nabla\Lambda(t;\alpha_0)^{\top}(\hat\alpha-\alpha_0)\\
 &\quad +n(\hat\alpha-\alpha_0)^{\top}\nabla\Lambda(t;\alpha_0)(\hat\alpha-\alpha_0)^{\top}\nabla^2\Lambda(t;\tilde\alpha)(\hat\alpha-\alpha_0)\\
 &\quad +\frac{n}{4}\left((\hat\alpha-\alpha_0)^{\top}\nabla^2\Lambda(t;\tilde\alpha)(\hat\alpha-\alpha_0)\right)^2.
\end{align*}
Using the consistency of $\hat\alpha-\alpha_0$ and the asymptotic normality of $\sqrt n(\hat\alpha-\alpha_0)$ from standard maximum likelihood theory, each of the last two terms in the above equation tends to $0$ in probability as $n$ tends to infinity. Therefore
\begin{align*}
& n\big(\Lambda(t;\hat\alpha)-\Lambda(t;\alpha_0)\big)^2-(n-1)\big(\Lambda(t;\hat\alpha^{(-l)})-\Lambda(t;\alpha_0)\big)^2\\
&\quad=(\hat\alpha^{(-l)}-\alpha_0)^{\top}\nabla\Lambda(t;\alpha_0)\nabla\Lambda(t;\alpha_0)^{\top}(n(\hat\alpha-\alpha_0)-(n-1)(\hat\alpha^{(-l)}-\alpha_0))\\
&\qquad +n(\hat\alpha-\hat\alpha^{(-l)})^{\top}\nabla\Lambda(t;\alpha_0)\nabla\Lambda(t;\alpha_0)^{\top}(\hat\alpha-\alpha_0)+o_{\mathbb P}(1)\\
&\quad=(\hat\alpha^{(-l)}-\alpha_0)^{\top}\nabla\Lambda(t;\alpha_0)\nabla\Lambda(t;\alpha_0)^{\top}(I^{-1}\nabla \log f(X_l;\alpha_0)+R_n)\\
&\qquad +(\alpha_0-\hat\alpha^{(-l)}+I^{-1}\nabla \log f(X_l;\alpha_0)+R_n')^{\top}\nabla\Lambda(t;\alpha_0)\nabla\Lambda(t;\alpha_0)^{\top}(\hat\alpha-\alpha_0)+o_{\mathbb P}(1),
\end{align*}
where the last two lines were derived from Equation~\eqref{eq:alpha} and $R_n$, $R_n'$ both tend to $0$ in probability. The consistency of $\hat\alpha$ and $\hat\alpha^{(-l)}$ shows that $A_{n,2}=o_{\mathbb P}(1)$. We now study the term $A_{n,1}$. From Equation~\eqref{eq:Lambda},
\begin{align*}
&n\big(\Lambda(t;\hat\alpha)-\Lambda(t;\alpha_0)\big)-(n-1)\big(\Lambda(t;\hat\alpha^{(-l)})-\Lambda(t;\alpha_0)\big)\\
&\quad=(n(\hat\alpha-\alpha_0)-(n-1)(\hat\alpha^{(-l)}-\alpha_0))^{\top}\nabla\Lambda(t;\alpha_0)\\
&\qquad + \frac 12 (\hat\alpha^{(-l)}-\alpha_0)^{\top}\nabla^2\Lambda(t;\tilde\alpha)(n(\hat\alpha-\alpha_0)-(n-1)(\hat\alpha^{(-l)}-\alpha_0))\\
&\qquad + \frac 12n(\hat\alpha-\hat\alpha^{(-l)})^{\top}\nabla^2\Lambda(t;\tilde\alpha)(\hat\alpha-\alpha_0).
\end{align*}
Using similar arguments as before, the last two lines of this Equation tend to $0$ in probability from Equation~\eqref{eq:alpha} and from the consistency of $\hat\alpha$ and $\hat\alpha^{(-l)}$. Finally, using again Equation~\eqref{eq:alpha}
\begin{align*}
A_{n,1}=-S(t;\alpha_0)\nabla \log f(X_l;\alpha_0)^{\top}I^{-1}\nabla\Lambda(t;\alpha_0)+o_{\mathbb P}(1).
\end{align*}
This equality combined with Equation~\eqref{eq:S} give
\begin{align*}
n S(t;\hat\alpha) - (n-1) S(t;\hat\alpha^{(-l)}) &=S(t;\alpha_0) - S(t;\alpha_0)\nabla \Lambda (t;\alpha_0)^{\top}{I}^{-1} \nabla \log f(X_l;\alpha_0)+o_{\mathbb P}(1).
\end{align*}
The final result of Proposition~\ref{paramsurv} is obtained by simply replacing each quantity by its consistent estimator. Integrating the equation in Proposition~\ref{paramsurv} directly yields the approximation for the RMST. By careful examination of the remainder term, we directly see that its integral over $[0,\tau]$ is also $o_{\mathbb P}(1)$.

\subsection{Log-likelihood, score vector and Hessian matrix in the piecewise constant hazard model}\label{sec:scoreHessian_pch}

In this section, we study the parametric piecewise constant hazard model defined as follows: 
$\lambda(t;\alpha)=\sum_{k=1}^K \alpha_k I_k(t)$ where $I_k(t)=I(c_{k-1}<t\leq c_k)$, $c_0=0<c_1<\cdots <c_K=+\infty$. The cumulative hazard function is then equal to
\begin{align*}
\Lambda(t;\alpha)=\sum_{k=1}^K \alpha_k (c_k\wedge t-c_{k-1}) I(c_{k-1}\leq t).
\end{align*} 
Under the mixed-case of interval-censored and exact data, we can directly write the log-likelihood as the sum between the log-likelihood of strictly interval-censored observations and the log-likelihood of exact observations. For the latter part see~\cite{aalen_borgan_book}. Recall that $X_i=(L_i,R_i)$ and $f(X_i;\alpha)$ denotes the density of the observations with parameter $\alpha$ evaluated at $X_i$. 
For strictly interval-censored data ($L_i\neq R_i$), the log-likelihood $\ell(\alpha)$ can be written as (see~\cite{sun07} or~\cite{bouaziz21})
\begin{align*}
\ell(\alpha)=\sum_{i=1}^n \log f(X_i;\alpha)&=\sum_{i=1}^n\bigg\{-(1-\Delta_i)\Lambda(L_i;\alpha)\\
&\qquad +\Delta_i\left(\log\Big(1-\exp\big(\Lambda(L_i;\alpha)-\Lambda(R_i;\alpha)\big)\Big)-\Lambda(L_i;\alpha)\right)\bigg\},
\end{align*}
where we used the notation $\Delta_i=I(R_i<+\infty)$ to denote uncensored observations. The $k^{\text{th}}$ component of the score vector is equal to:
\begin{align}\label{eq:score}
\frac{\partial \ell(\alpha)}{\partial\alpha_k}=\sum_{i=1}^n \frac{\partial \log f(X_i;\alpha)}{\partial \alpha_k}& = \sum_{i=1}^n\Big\{-(c_k\wedge L_i-c_{k-1})I(c_{k-1}\leq L_i) \nonumber \\
&\qquad +\Delta_i\frac{(c_k\wedge R_i-c_{k-1})I(c_{k-1}\leq R_i)-(c_k\wedge L_i-c_{k-1})I(c_{k-1}\leq L_i)}{1-\exp\big(\Lambda(L_i;\alpha)-\Lambda(R_i;\alpha)\big)}\nonumber\\
&\qquad \quad\times\exp\big(\Lambda(L_i;\alpha)-\Lambda(R_i;\alpha)\big)\Big\}\cdot
\end{align}
The $k\times k'$ component of the Hessian matrix is equal to:
\begin{align}\label{eq:Hessian}
\frac{\partial^2 \ell(\alpha)}{\partial\alpha_{k'}\partial\alpha_k}& = -\sum_{i=1}^n\Delta_i\left\{\frac{(c_k\wedge R_i-c_{k-1})I(c_{k-1}\leq R_i)-(c_k\wedge L_i-c_{k-1})I(c_{k-1}\leq L_i)}{1-\exp\big(\Lambda(L_i;\alpha)-\Lambda(R_i;\alpha)\big)}\right.\nonumber\\
&\qquad \qquad\times\big((c_k'\wedge R_i-c_{k'-1})I(c_{k'-1}\leq R_i)-(c_k'\wedge L_i-c_{k'-1})I(c_{k'-1}\leq L_i)\big)\nonumber\\
&\qquad \qquad\times\exp\big(\Lambda(L_i;\alpha)-\Lambda(R_i;\alpha)\big)\nonumber\\
&\qquad \quad +\frac{(c_k\wedge R_i-c_{k-1})I(c_{k-1}\leq R_i)-(c_k\wedge L_i-c_{k-1})I(c_{k-1}\leq L_i)}{\Big(1-\exp\big(\Lambda(L_i;\alpha)-\Lambda(R_i;\alpha)\big)\Big)^2}\nonumber\\
&\qquad \qquad\times \big((c_k'\wedge R_i-c_{k'-1})I(c_{k'-1}\leq R_i)-(c_k'\wedge L_i-c_{k'-1})I(c_{k'-1}\leq L_i)\big)\nonumber\\
&\qquad \qquad\times\exp\Big(2\big(\Lambda(L_i;\alpha)-\Lambda(R_i;\alpha)\big)\Big)\Bigg\}\cdot
\end{align}
The Fisher information is equal to the expectation of minus the Hessian matrix divided by $n$. Looking at its expression, we directly see that a necessary condition for the Fisher information to be positive definite is to assume that
\begin{align*}
\mathbb P(\Delta=1, [L,R] \cap (c_{k-1},c_k]\neq \emptyset) >0, \forall k=1,\ldots,K.
\end{align*}
Another important condition for the model to be identifiable is to assume that $\mathbb E_{\alpha_0}[f(X;\alpha)]$ has a unique maximum with respect to $\alpha$, equal to $\alpha_0$, where the notation $\mathbb E_{\alpha_0}$ means the expectation is taken with respect to the true parameter $\alpha_0$. However, it is clear from Equation~\eqref{eq:score} that $\mathbb E_{\alpha_0}[\partial f(X;\alpha)/\partial \alpha]$ cannot vanish if $\mathbb P(L>c_{k-1})=0$. Therefore, a second necessary condition for the model to be identifiable is to assume
\begin{align*}
\mathbb P(L>c_{k-1}) >0, \forall k=1,\ldots,K.
\end{align*}
Those two conditions have opposite effects on the estimation method if they are violated. In case the first one is not valid for a given $k$ then it will not be possible to compute the corresponding estimator $\hat\alpha_k$ from the Newton-Raphson algorithm (since the Hessian will not be invertible) while using the EM algorithm (which does not involve the Score vector nor the Hessian matrix), the estimator $\hat\alpha_k$ will become smaller at each iteration step until eventually reaching the value $0$. This situation can be numerically resolved in the latter case, by simply setting the iterated estimate $\hat\alpha_k$ to $0$ when it reaches a value below a fixed threshold. However this situation is problematic for the computation of the pseudo-values. This can be easily seen by recalling that pseudo-values should average to the initial estimator. In Proposition~\ref{paramsurv} this simply follows from the fact that $\sum_{l=1}^n \nabla \log f(X_l;\hat\alpha)=0$ from regularity conditions for maximum likelihood estimation. However, the $k^{\text{th}}$ component of the score vector will never vanish if the first condition is not valid, leading to incorrect pseudo-values. 

On the other hand, if the second condition is not valid for a given $k$, the algorithm will attempt to minimise the term $\exp(-\Lambda(R_i;\alpha))$ from Equation~\eqref{eq:score} and as a consequence the corresponding estimator $\hat\alpha_k$ will become larger at each iteration step of the EM algorithm, diverging to infinity.

Finally, note that if the log-likelihood only include exact observations $L=R$, the conditions then translate to $\mathbb P(c_{k-1}<L<c_k) >0, \forall k=1,\ldots,K.$

\subsection{Implementation of the pseudo-observations for the survival function and the RMST in the pch model}\label{sec:details_pseudoIC}

In this section we provide the precise expression of the terms involved in Proposition~\ref{paramsurv} for the pch model. We have 
\begin{align*}
S(t;\alpha)&=\exp\Big(-\sum_{k=1}^K \alpha_k(t\wedge c_k-c_{k-1})I(c_{k-1}\leq t)\Big)\\
\frac{\partial\Lambda(t;\alpha)}{\partial\alpha_k}&=(c_k\wedge t-c_{k-1}) I(c_{k-1}\leq t),
\end{align*} 
while the expression of the gradient of the density $\nabla \log f(X_l;\alpha)$ is given by the term between brackets in Equation~\eqref{eq:score} and $\hat I$ is equal to minus the Hessian matrix (see Equation~\eqref{eq:Hessian}) divided by $n$.

For the integrated version we need to precise how to compute the integral between $0$ and $\tau$ of $S(t;\alpha)$ and the integral between $0$ and $\tau$ of $S(t;\alpha)\nabla \Lambda(t;\alpha)$. We first notice that
\begin{align*}
S(t;\alpha)&=\exp\Big(-\sum_{k=1}^K \alpha_k(c_k-c_{k-1})I(c_{k}\leq t)\Big)\exp\Big(-\sum_{k=1}^K \alpha_k(t-c_{k-1})I(c_{k-1}\leq t\leq c_k)\Big),
\end{align*}
and
\begin{align*}
\int_0^{\tau} S(t;\alpha)dt &= \sum_{l=1}^K \int_{c_{l-1}}^{c_l\wedge \tau}\!\! S(t;\alpha)dt\; I(\tau>c_{l-1})\\
&= \sum_{l=1}^K \int_{c_{l-1}}^{c_l\wedge \tau}\!\!\exp\Big(-\sum_{k=1}^K \alpha_k(c_k-c_{k-1})I(c_{k}\leq t)\Big)\exp\Big(- \alpha_l(t-c_{l-1})\Big)dt\,I(\tau>c_{l-1}).
\end{align*}
Set $A_1=0$ and for $l\geq 2$, define
\begin{align*}
A_l=-\sum_{k=1}^{l-1} \alpha_k(c_k-c_{k-1})+c_{l-1}\alpha_l.
\end{align*}
For the first term we now have:
\begin{align*}
\int_0^{\tau} S(t;\alpha)dt &=\sum_{l=1}^K \exp(A_l)\int_{c_{l-1}}^{c_l\wedge \tau}\exp(- \alpha_lt)dt\,I(\tau>c_{l-1})\\
&=\sum_{l=1}^K \exp(A_l)\alpha_l^{-1}\Big(\exp(- \alpha_lc_{l-1})-\exp\big(-\alpha_l (c_l\wedge \tau)\big)\Big)\,I(\tau>c_{l-1}).
\end{align*}
For the second term we have:
\begin{align*}
\int_0^{\tau} S(t;\alpha)\frac{\partial\Lambda(t;\alpha)}{\partial\alpha_k}dt &=\sum_{l=1}^K \int_{c_{l-1}}^{c_l\wedge \tau}\!\!S(t;\alpha)(c_k\wedge t-c_{k-1}) I(c_{k-1}\leq t)dt\,I(\tau>c_{l-1})\\
&=\sum_{l=k+1}^{K} \int_{c_{l-1}}^{c_{l}\wedge \tau}\!\!S(t;\alpha)dt\,(c_k-c_{k-1})I(\tau>c_{l-1})\\
&\quad +\int_{c_{k-1}}^{c_{k}\wedge \tau}\!\!tS(t;\alpha)dt\,I(\tau>c_{k-1})- c_{k-1}\int_{c_{k-1}}^{c_{k}\wedge \tau}\!\!S(t;\alpha)dt\,I(\tau>c_{k-1}).
\end{align*}
From the previous calculation on the first term, we easily see that on the one hand
\begin{align*}
\int_{c_{l-1}}^{c_{l}\wedge \tau}\!\!S(t;\alpha)dtI(\tau>c_{l-1}) = \exp(A_l) \alpha_l^{-1}\Big(\exp(- \alpha_lc_{l-1})-\exp\big(-\alpha_l (c_l\wedge \tau)\big)\Big)\,I(\tau>c_{l-1}).
\end{align*}
%
On the other hand, we have:
\begin{align*}
& \int_{c_{l-1}}^{c_{l}\wedge \tau}\!\!tS(t;\alpha)dtI(\tau>c_{l-1})\\
 & = \exp(A_l)\int_{c_{l-1}}^{c_{l}\wedge \tau}\!\!t\exp(-t\alpha_l)dtI(\tau>c_{l-1})\\
& = \exp(A_l)\bigg(\alpha_l^{-2}\Big(\exp\big(-c_{l-1}\alpha_l\big)-\exp\big(-(c_{l}\wedge \tau)\alpha_l\big)\Big)\\
&\quad+\alpha_l^{-1}\Big(c_{l-1}\exp\big(-c_{l-1}\alpha_l\big)-(c_{l}\wedge \tau)\exp\big(-(c_{l}\wedge \tau)\alpha_l\big)\Big)\bigg)I(\tau>c_{l-1}),
\end{align*}
where the last equation was obtained using integration by parts. Gathering all elements allows to implement the second equation in Proposition~\ref{paramsurv}.

\subsection{Proof of Proposition~\ref{paramsurv_approx_property}}

First, from a Taylor expansion of $\int_0^{\tau} S(t;\alpha_z)dt$ around $\int_0^{\tau} S(t;\alpha_0)dt$, we obtain:
\begin{align*}
\mathbb E(T^*_l\wedge \tau\mid Z_l=z) & = \int_0^{\tau} S(t;\alpha_z)dt\\
 & = \int_0^{\tau} S(t;\alpha_0)dt +  \int_0^{\tau} (\nabla S(t;\alpha_0))^{\top}dt(\alpha_z-\alpha_0)\\
 &\quad + R_{1,z},
\end{align*}
with
\begin{align*}
R_{1,z}= \frac 12 (\alpha_z-\alpha_0)^{\top} \int_0^{\tau} \nabla^2 S(t;\tilde\alpha_z)dt(\alpha_z-\alpha_0)
\end{align*}
and $\tilde\alpha_z$ is on the real line between $\alpha_0$ and $\alpha_z$. Then, since $\alpha_z$ maximises with respect to $\alpha$ the expected log-likelihood $\mathbb E(\log f(X_l;\alpha)\mid Z_l=z))$ we have that $\mathbb E(\nabla\log f(X_l;\alpha_z)\mid Z_l=z))=0$. Then, from a Taylor expansion around $\alpha_0$ we obtain:
\begin{align*}
0=\mathbb E(\nabla\log f(X_l;\alpha_0)\mid Z_l=z))+(\alpha_z-\alpha_0)^{\top}\mathbb E(\nabla^2 \log f(X_l;\tilde{\tilde\alpha}_z)\mid Z_l=z)),
\end{align*}
where $\tilde{\tilde\alpha}_z$ is on the real line between $\alpha_0$ and $\alpha_z$. As a result, we have:
\begin{align*}
\alpha_z-\alpha_0&=I^{-1}\mathbb E(\nabla \log f(X_l;\alpha_0)\mid Z_l=z)\\
&\quad+\left\{(-\mathbb E(\nabla^2 \log f(X_l;\tilde{\tilde\alpha}_z)\mid Z_l=z)))^{-1}-I^{-1}\right\}\mathbb E(\nabla \log f(X_l;\alpha_0)\mid Z_l=z).
\end{align*}
Gathering all parts, we have proved that
\begin{align*}
\mathbb E(T^*_l\wedge \tau\mid Z_l=z) &= \int_0^{\tau} S(t;\alpha_0)dt + \int_0^{\tau} (\nabla S(t;\alpha_0))^{\top}dt I^{-1}\mathbb E(\nabla \log f(X_l;\alpha_0)\mid Z_l=z)\\
&\quad+R_{1,z}+R_{2,z},
\end{align*}
with
\begin{align*}
R_{2,z}=\int_0^{\tau} (\nabla S(t;\alpha_0))^{\top}dt\left\{(-\mathbb E(\nabla^2 \log f(X_l;\tilde{\tilde\alpha}_z)\mid Z_l=z)))^{-1}-I^{-1}\right\}\mathbb E(\nabla \log f(X_l;\alpha_0)\mid Z_l=z).
\end{align*}
Finally, writing $S(t;\alpha_0)=\exp(-\Lambda(t;\alpha_0))$, we directly obtain
\begin{align*}
\int_0^{\tau} \nabla S(t;\alpha_0)dt=-\int _0^{\tau}\!S(t;\alpha_0)\nabla \Lambda (t;\alpha_0)dt,
\end{align*}
which concludes the proof.

\subsection{Proof of Proposition~\ref{paramsurv_pch_property}}

Let $X_i=(T_i,\Delta_i)$ and assume the pch model for $\lambda$. In this model, the cumulative hazard function evaluated at $T_i$, is equal to
\begin{align*}
\Lambda(T_i;\alpha)=\sum_{k=1}^K \alpha_k (c_k\wedge T_i-c_{k-1}) I(c_{k-1}\leq T_i).
\end{align*}
However, since we want to prove the result when the mesh of the partition $0=c_0<c_1<\cdots<c_K=+\infty$ tends to zero, we can write without loss of generality that there exists a $\delta>0$, such that for any partition whose mesh is less than $\delta$ we have for all $T_i$, $i=1,\ldots,n$,
\begin{align*}
\Lambda(T_i;\alpha)=\sum_{k=1}^K \alpha_k (c_k-c_{k-1}) I(c_{k}\leq T_i).
\end{align*}
Therefore, for a partition $0=c_0<c_1<\cdots<c_K=+\infty$ such that $\max_k|c_k-c_{k-1}|<\delta$,
\begin{align*}
\log f(X_i;\alpha) & = \Delta_i\sum_{k=1}^K\log(\alpha_k) I_k(T_i)-\sum_{k=1}^K \alpha_k (c_k-c_{k-1}) I(c_{k}\leq T_i),\\
\frac{\partial}{\partial\alpha_k}\log f(X_i;\alpha)& =\frac{\Delta_i}{\alpha_k}I_k(T_i)-(c_k-c_{k-1})I(c_{k}\leq T_i),\\
\frac{\partial^2}{\partial\alpha_k^2}\log f(X_i;\alpha)& =-\frac{\Delta_i}{\alpha_k^2}I_k(T_i),\\
\frac{\partial^2}{\partial\alpha_k'\partial \alpha_k}\log f(X_i;\alpha)& =0, \text{ for }k\neq k'.
\end{align*}
We therefore have that the Fisher information $I$ is a diagonal matrix whose $k$th element, $k=1,\ldots,K$, is equal to $\mathbb E(\Delta_iI_k(T_i))/(\alpha^0_k)^2$, where $\alpha^0_k$ is the $k$th component of $\alpha_0$. Also under standard maximum likelihood regularity conditions, the true parameter $\alpha_0$ verifies that
\begin{align*}
\mathbb E\left(\frac{\partial}{\partial\alpha_k} \log f(X_i;\alpha_k^0)\right)=0,
\end{align*}
which is equivalent to
\begin{align*}
\alpha_k^0 = \frac{\mathbb E(\Delta_i I_k(T_i))}{(c_k-c_{k-1})\mathbb P(T_i\geq c_{k} )}\cdot
\end{align*}
Now, we have that $\mathbb E(\Delta_lI_k(T_l))/(c_k-c_{k-1})$ tends to $H'_1(c_{k-1})$ as $(c_k-c_{k-1})$ tends to $0$, where $H_1(t)=\mathbb P(T\leq t,\Delta=1)$. Therefore, as the limit $(c_k-c_{k-1})$ goes to $0$, $\alpha_k^0$ tends to $\lambda(c_k)$, the true hazard rate evaluated at $c_k$. In other words, the parametric hazard $\lambda(t;\alpha_0)$ tends to the true hazard function $\lambda(t)$ as $(c_k-c_{k-1})$ goes to $0$.

Then, $\nabla\Lambda(t;\alpha_0)$ is a vector whose $k$th component is equal to $(c_k-c_{k-1})I(c_k\leq t)$. As a consequence,
\begin{align*}
& \nabla \Lambda (t;\alpha_0)^{\top}{ I}^{-1} \nabla \log f(X_l;\alpha_0)\\
& = \sum_{k=1}^K I(c_k\leq t)\frac{(\alpha_k^0)^2(c_k-c_{k-1})}{\mathbb E(\Delta_lI_k(T_l))}\left(\frac{\Delta_l}{\alpha_k^0}I_k(T_l)-(c_k-c_{k-1})I(c_k\leq T_l)\right).
\end{align*}
 In the first term of this equation, there can only be one interval of the form $[c_{k-1},c_k]$ that contains $T_l$ and therefore the sum over $k$ is not zero for only one value of $k$. As the limit of $(c_k-c_{k-1})$ goes to $0$, $T_l$ gets closer to $c_k$ and $c_{k-1}$ and $\alpha_k^0$ tends to $\lambda(T_l)$. More precisely,
\begin{align*}
\lim_{c_k-c_{k-1}\to 0}\sum_{k=1}^K I(c_k\leq t)\frac{(\alpha_k^0)^2(c_k-c_{k-1})}{\mathbb E(\Delta_lI_k(T_l))}\frac{\Delta_l}{\alpha_k^0}I_k(T_l)   = \frac{\Delta_l I(T_l\leq t) \lambda(T_l)}{H'_1(T_l)}\cdot
\end{align*}
On the other hand, the second term of the equation is simply a Riemann integral. We have:
\begin{align*}
\lim_{c_k-c_{k-1}\to 0}\sum_{k=1}^K I(c_k\leq t)\frac{(\alpha_k^0)^2(c_k-c_{k-1})^2}{\mathbb E(\Delta_lI_k(T_l))}I(c_k\leq T_l) = \int_0^{T_l\wedge t} \frac{(\lambda(u))^2}{H'_1(u)}du\cdot
\end{align*}
Noticing that $\lambda(t)=H_1'(t)/H(t)$, with $H(t)=\mathbb P(T\geq t )$, we therefore have proved
\begin{align*}
\lim_{c_k-c_{k-1}\to 0}\nabla \Lambda (t;\alpha_0)^{\top}{ I}^{-1} \nabla \log f(X_l;\alpha_0)= \frac{\Delta_l I(T_l\leq t)}{H(T_l)}-\int_0^{T_l\wedge t} \frac{dH_1(u)}{(H(u))^2}\cdot
\end{align*}
We now want to compute the conditional expectation with respect to $Z_l$ of this quantity. Let $H^Z_1(t)=\mathbb P(T\leq t,\Delta=1\mid Z)$, $H^Z(t)=\mathbb P(T\geq t \mid Z)$ and write
\begin{align*}
\mathbb E\left(\frac{\Delta_l I(T_l\leq t)}{H(T_l)}-\int_0^{T_l\wedge t} \frac{dH_1(u)}{(H(u))^2}\mid Z_l\right)&=\int_0^t \frac{dH_1^{Z_l}(u)}{H(u)}-\int_0^{t} \frac{H^{Z_l}(u)dH_1(u)}{(H(u))^2}\\
&=\int_0^t \left(\frac{(H_1^{Z_l}(u))'}{H^{Z_l}(u)}-\frac{(H_1(u))'}{H(u)}\right)\frac{H^{Z_l}(u)}{H(u)}du\\
&=\int_0^t (\lambda^{Z_l}(u)-\lambda(u))\frac{H^{Z_l}(u)}{H(u)}du,
\end{align*} 
where $\lambda^{Z_l}$ represents the conditional hazard function. Setting $q_{Z_l}(u)=(\lambda^{Z_l}(u)-\lambda(u))H^{Z_l}(u)/H(u)$ as in~\cite{jacobsen2016note} we obtain from their proofs in Section A.2, Equation (11), that $-S^{Z_l}/S$ is a primitive of $q_{Z_l}$, where $S^{Z_l}$ is the conditional survival function of $T^*$. As a consequence,
\begin{align*}
\int_0^t (\lambda^{Z_l}(u)-\lambda(u))\frac{dH^{Z_l}(u)}{H(u)}=1-\frac{S^{Z_l}(t)}{S(t)}\cdot
\end{align*}
Gathering all the parts, we have proved that
\begin{align*}
\lim_{c_k-c_{k-1}\to 0}\mathbb E\left(\varphi(X_l;\alpha_0)\mid Z_l\right)&=\int _0^{\tau}\!S(t)dt - \int _0^{\tau}\!S(t)\left(1-\frac{S^{Z_l}(t)}{S(t)}\right)dt\\
&=\int _0^{\tau}\! S^{Z_l}(t) dt = \mathbb E(T_l\wedge \tau\mid Z_l),
\end{align*}
which concludes the proof.


\bibliographystyle{unsrt}
\bibliography{biblio}

\end{document}